\documentclass[preprint,12pt]{elsarticle}
\usepackage[top=1in, bottom=1.25in, left=1.25in, right=1.25in]{geometry}

\usepackage{amssymb}
\usepackage{amsmath}
\usepackage{amsthm}

\usepackage{ART-SFVEM}
\input{Pack.sty}
\newtheorem{problem}{Problem}
\newtheorem{remark}{Remark}
\newtheorem{lemma}{Lemma}
\newtheorem{proposition}{Proposition}
\newtheorem{theorem}{Theorem}

\newcommand{\mydefeq}{\vcentcolon=}

\begin{document}

\begin{frontmatter}



\title{Stabilization-Free General Order Virtual Element Methods for Neumann Boundary Optimal Control Problems in Saddle Point Formulation}


\author[poli,gncs]{Andrea Borio} 
\author[poli,gncs]{Francesca Marcon}
\author[poli,gncs]{Maria Strazzullo}


\affiliation[poli]{organization={Politecnico of Torino, Department of Mathematical Sciences ``Giuseppe Luigi Lagrange''},
            addressline={Corso Duca degli Abruzzi, 24}, 
            city={Torino},
            postcode={10129}, 
            state={},
            country={Italy}}
\affiliation[gncs]{organization={INdAM research group GNCS member},
  country={Italy}}

\begin{abstract}
  In this work, we explore the application of Stabilization-Free Virtual Element Methods for Neumann
  boundary Optimal Control Problems in saddle point formulation. The method is proposed for
  arbitrary polynomial order of accuracy and general polygonal meshes. Our contribution includes a
  rigorous \emph{a priori} error estimate that holds for general polynomial order. On the numerical
  side, we present (i) an initial convergence test that reflects our theoretical findings, (ii) a
  second test analyzing the role of the stabilization term in the Virtual Element Method (VEM)
  formulation and its influence on the approximation error, and (iii) a third test case based on a
  more application-oriented experiment.  The stabilization-free approach is proposed as an
  alternative strategy to circumvent issues related to the choice of the stabilization parameter in
  standard VEM formulations.
\end{abstract}


\end{frontmatter}


\section{Introduction}
\label{sec:intro}

In this paper, we consider a Stabilization-Free Virtual Element Method (SFVEM) approximation for optimal Neumann boundary
control. Optimal Control Problems (OCPs) governed by partial differential
equations (PDEs) are widespread in many practical applications in scientific and
engineering contexts, see e.g. \cite{Ballarin2022307,Pichi20221361,Strazzullo2023,Strazzullo2021841} and the references therein. The main goal
of OCPs is to reach a beneficial behavior of the PDE solution by employing an
external variable, the \emph{control}. The role of the control consists in
making the \emph{state} solution similar to a desired configuration employing a functional minimization. The control is
\emph{distributed} when its action is defined all over the domain, while we talk
about \emph{boundary control} if the control affects the boundary
conditions. Due to its coupled nature, the numerical analysis for OCPs
governed by PDEs is a fascinating mathematical challenge. Several works
have been proposed with Finite Element Methods, Discontinuous Galerkin Methods,
Finite Volumes and hybrid methods (see for instance \cite{casas2008error,kumar2017discontinuous,YUCEL20152414}).

Lately, great attention has been paid to VEM, a family
of polygonal methods for simulating PDEs. VEM
can be applied to arbitrary polygonal meshes, guaranteeing huge flexibility in
terms of geometrical complexity, and they might play a crucial role in many
real-life scenarios.  Indeed, after being conceived in the seminal paper
\cite{Beirao2013a}, they have been applied to various contexts of interest.
The investigation of boundary OCPs through VEM
is still limited, yet of utmost importance. Indeed, the versatility of arbitrary
polygonal meshes, combined with the action of control strategies on complex
geometries, might not only be a challenging theoretical topic but also relevant
from a practical viewpoint.  To our knowledge, $C^1$-VEM is presented in
\cite{BrennerC1} for distributed OCPs governed by elliptic equations and with
a constrained state variable. The VEM discretization for distributed control is
the topic of several investigations, see e.g.
\cite{Sun20242019,Tushar2022,TUSHAR2023134,Wang2021}.  Some VEM
works related to boundary OCPs are \cite{Tushar2024_diri} for the Dirichlet
case with control constraints and \cite{Liu2023} for the Neumann case with
control constraints and approximation of the first order.
A key ingredient of the standard VEM approximation for any problem of interest is the use of discrete bilinear forms composed of a consistent polynomial part and a stabilization term ensuring coercivity.
It is known that the choice of stabilization remains somewhat arbitrary due to its problem dependance, and its influence on the system cannot be reliably controlled in more complex settings, such as the coupled problems and control applications considered here.
 In the literature, many alternatives have been proposed that do not rely on the stabilization term.
 Starting from the works \cite{DAltri2021, BBME2VEM,berrone_comparison_2022}, a research line has grown, defining the so-called \textit{Stabilization-Free} VEM (SFVEM) and focusing on the definition of VEM bilinear forms that are self‑stabilized through the use of higher‑order polynomial projections (see for instance \cite{Lamperti2022,CHEN202388,BERRONE2023108641,XU2023,BLMVMixedStabFree}). For further studies on the role of stabilization and SFVEM, the reader may refer to \cite{russo2023quantitative,Beirao2023}.
 In this contribution, we investigate the SFVEM within the framework of OCPs, employing the strategy proposed in \cite{Berrone2024,BERRONE2025117839} which exploits the properties of divergence‑free polynomial spaces.
In particular, we focus on the unconstrained Neumann boundary control problem in a
\emph{saddle point} structure.  
This naturally
derives from the optimality conditions and considers the coupled nature of the
system leading to more robust and stable systems, solving for all the variables
at once.  
This formulation is generally more effective
than iterative algorithms. For these kinds of coupled systems, SFVEM might represent an asset in terms of robustness of the method compared to stabilized VEM. 
To the best of our knowledge, we
differentiate from previous contributions in the following terms:
\begin{itemize}
\item the presented {VEM scheme does not require any additional stabilization term and has arbitrary polynomial accuracy;}
\item we analyze the problem in saddle point formulation as presented in
  \cite{manzoniOCP}, {hence providing a complete proof of well-posedness and \emph{a priori} error estimates;}
\item the numerical tests cover high-order schemes, and we numerically investigate how the stabilization parameter affects the
  method's accuracy, {assessing the flexibility and robustness of the proposed stabilization-free scheme.}
\end{itemize}

The paper is structured as follows: in Section \ref{sec:problem} we introduce
the linear-quadratic Neumann boundary OCP at the continuous level in a saddle
point framework. Section \ref{sec:VEM} focuses on the Stabilization-Free Virtual Element Method (SFVEM) approximation of such
optimization problem and proves the well-posedness of the saddle point
structure. In Section \ref{sec:err_estimates} we derive optimal a priori error estimates.
{Section \ref{sec:results} presents numerical results that confirm the expected convergence behavior, we also analyze how the accuracy of the stabilized VEM formulation depends on the choice of the stabilization parameter, and we present the proposed SFVEM method as a robust strategy for a more application-oriented test.}

Finally, conclusions follow in Section
\ref{sec:conclusions}.

Throughout the work, $(\cdot,\cdot)_{\omega}$ denotes the standard
$\mathrm{L}^2$ scalar product defined on a generic $\omega\subset\mathbb{R}^2$
and inside the proofs, the symbol $C$ denotes any constant independent of the
mesh size.

\section{Problem Formulation}
\label{sec:problem}

Let us consider $\Omega \subset \mathbb R^2$ as a bounded domain with a polygonal boundary. We denote as $\Gamma_D$ and $\Gamma_N$ the portions of the boundary $\partial \Omega$ where Dirichlet and Neumann conditions apply. Moreover, we define a \emph{control boundary} $\Gamma_C$, where the control acts in form of Neumann condition. We remark that $\Gamma_D \cup \Gamma_N \cup \Gamma_C = \partial \Omega$ and they are pairwise disjoint. Moreover, for simplicity we assume $\Gamma_D\neq\emptyset$.

The state variable $y$ and the control variable $u$ are considered, respectively, in 
$$
Y = \{y \in H^1(\Omega) \; | \; y = 0 \text{ on } \Gamma_D\} \text{ and } U = L^2(\Gamma_C).
$$
The case of non-homogeneous boundary conditions can be tackled by standard lifting procedures. 
Given a desired state $\dState \in Y_{\text{d}} = L^2(\OmegaObs)$, with $\OmegaObs \subset \Omega$, a forcing term $f \in L^2(\Omega)$ 
and a penalization parameter $\alpha > 0$ on the control action, we solve the following problem.
\begin{problem}\label{P1:OCP}
Find $(y,u) \in Y \times U$ such that
\begin{equation*}
J(y,u) = \min_{(w,v) \in Y \times U}  J(w,v) = \min_{(w,v) \in Y \times U} 
\half \norm[\lebl{\OmegaObs}]{w-\dState}^2 +\frac{ \alpha}{2} \norm[\lebl{\Gamma_C}]{v}^2 \,,
\end{equation*}
subject to the following PDE
\begin{equation}
\label{eq:state_eq}
\begin{cases}
\nabla \cdot (-\K \nabla w + \beta w) + \gamma w = f & \text{in } {\Omega}, \\
\displaystyle w = 0 & \text{on  } {\Gamma_D}, \\
\displaystyle \dn{w} = v & \text{on  } {\Gamma_C},\vspace{0.5mm} \\
\displaystyle \dn{w} = 0 & \text{on  } {\Gamma_N}, \\
\end{cases}
\end{equation}
with $\beta \in (L^{\infty}(\Omega))^2$ is divergence-free and $\gamma \in L^{\infty}(\Omega)$ is non negative on $\Omega$. Moreover, we consider the symmetric diffusive tensor $\mathcal K \in (L^{\infty}(\Omega))^{2 \times 2}$ satisfying  
\begin{equation}
\mathcal{K}_0 \, |\nu|^2 \;\le\; \nu \cdot \mathcal{K}(x)\, \nu \;\le\; \mathcal{K}_1 \, |\nu|^2,
\quad \forall\, \nu \in \mathbb{R}^2,\ \forall\, x \in \Omega,    
\end{equation}
for some positive constants  $\mathcal{K}_0$ and $\mathcal{K}_1$ and $| \cdot | $ denoting the Euclidean norm.
\end{problem}

Problem~\ref{P1:OCP} can be expressed in weak formulation as
find the pair $(y,u) \in Y \times U$ which minimizes $J(y,u)$ and satisfies 
\begin{equation}
\label{eq:state_eq_weak}
{a (y, q)} =
\scal[\Gamma_C]{u}{q} + \scal[\Omega]{f}{q}
\quad \forall q \in Y,
\end{equation}
where 
$a \goesto{Y}{Y}{\mathbb R}$ is given by 
\begin{equation}
  \label{eq:def-a}
a(y,q) \mydefeq
\scal[\Omega]{\K\nabla y}{\nabla q} + \scal[\Omega]{\beta \nabla y}{q}+ \scal[\Omega]{\gamma y}{q}
\, \quad \forall y,q \in Y.
\end{equation}

Since $a(\cdot, \cdot)$ is continuous and coercive and $\scal[\Gamma_C]{u}{q}$ is continuous,  \eqref{eq:state_eq_weak} is well-posed. 
Let $\ennorm[Y]{\cdot}$ denote the following norm: for each $w\in Y$
\begin{equation}
\ennorm[Y]{w}^2\eqdot\norm[\lebl{\Omega}]{\sqrt{\K}\nabla w}^2. 
\end{equation} 
\begin{remark}
\label{rem:H1norms}
For any $w\in Y$, there exists $\tilde C>0$ such that
\begin{equation}
    \norm[\sobh{1}{\Omega}]{w} \coloneqq \norm[\lebl{\Omega}]{\nabla w}^2 + \norm[\lebl{\Omega}]{w}^2 \leq \tilde C \ennorm[Y]{w}
\end{equation} 
holds true, due to a Poincaré inequality and the regularity hypothesis of $\mathcal K$. Thanks to the same hypothesis, we have that, for a positive constant $C$
\begin{equation}
    \norm[\lebl{\Omega}]{w}^2 \leq C \ennorm[Y]{w}.
\end{equation} 
holds true.
\end{remark}

We define the \emph{weak minimization problem} as 
\begin{equation}
\label{eq:min_weak}
(y,u)=\argmin\{J(w, v)  :  (w,  v) \in Y \times U \text{ and \eqref{eq:state_eq_weak} is verified}\}.
\end{equation}
Considering $\alpha>0$ and thanks to the assumptions over the bilinear forms $a(\cdot, \cdot)$ and $\scal[\Gamma_C]{u}{q}$, the minimization problem is well-posed and  a Lagrangian approach can be exploited to solve it \cite{manzoniOCP}.
First, we define an \emph{adjoint variable} $p \in Y$ and the Lagrangian functional as
\begin{align}
\label{eq:LG}
 \Lg (y,u,p) = J(y, u) + \displaystyle a (y, p)- \scal[\Gamma_C]{u}{p} -  \scal[\Omega]{f}{p}.
\end{align}
The pair $(y,u)$ minimizing \eqref{eq:min_weak}, is obtained solving the \emph{optimality system} that finds $(y, u, p) \in Y \times  U \times Y$ such that
\begin{equation}
\label{eq:optimality_system}
\begin{cases}
D_y\Lg(y, u, p)[w] = 0 & \forall w \in Y,\\
D_u\Lg(y, u, p)[v] = 0 & \forall v \in U,\\
D_p\Lg(y, u, p) [q]= 0 & \forall q \in Y,\\
\end{cases}
\end{equation}
where $D_{y}, D_{u}$ and $D_{p}$ denotes the differentiation with respect to the three variables, respectively.
In this way, the PDE-constrained optimization problem is reframed in an unconstrained optimization and the minimum of \eqref{eq:min_weak} is represented by the state and control components of the stationary point of \eqref{eq:LG}. The strong formulation of the optimality system reads: 
\begin{problem}\label{PB:OS}
Find $(y,u,p) \in Y \times U \times Y$ such that
\begin{equation}
\label{eq:strong_os}
\begin{cases}
\displaystyle y \ChObs + \nabla \cdot (-\K \nabla p-\beta p) +\gamma p =
\dState\ChObs & \text{ in } \Omega, \quad \text{\textit{adjoint equation}}\\
\alpha u -  p = 0 & \text{ in } \Gamma_C, \;\; \text{\textit{optimality condition}}\\
\nabla \cdot (-\K\nabla y+\beta y)  + \gamma y  = f & \text{ in } \Omega, \quad \text{\textit{state constraint}}\\
\displaystyle \frac{\partial y}{\partial n} = u \quad \text{and} \quad \frac{\partial p}{\partial n} + \beta n p= 0 & \text{ on }\Gamma_C, \vspace{2mm}\\ 
\displaystyle \dn{y} = 0 \quad \text{and} \quad \frac{\partial p}{\partial n} + \beta n p = 0 &\text{ on } \Gamma_N,\\
y = p =0 &  \text{ on } \Gamma_D, \\
\end{cases}
\end{equation}
where $\ChObs$ is the indicator function of the observation domain $\OmegaObs$. 
\end{problem}

Let us consider $\bs{x}:=(y,u)\in X$, where $X$ is the product space between $Y$ and $U$ equipped with the scalar product 
\begin{equation}
\scal[X]{\bs{x}_1}{\bs{x}_2}=\scal[\Omega]{y_1}{y_2}+\scal[\Gamma_C]{u_1}{u_2}\,,
\end{equation}
for any $\bs{x}_1=(y_1,u_1)$ and $\bs{x}_2=(y_2,u_2)$ in $X$,
and with the norm
\begin{equation}
\norm[X]{\bs{x}}^2 = \ennorm[Y]{y}^2 + \norm[\lebl{\Gamma_C}]{u}^2,
\end{equation}
for any $\bs{x}=(y,u)$.

Then, introducing the bilinear forms
$\A{}{}:X\times X\to\mathbb{R}$
and
$\B{}{}:Y\times X\to\mathbb{R}$
given by
\begin{equation}
\A{\bs{x}}{\bs{\varsigma}} = \scal[{\OmegaObs}]{y}{w}  + \alpha \scal[\Gamma_C]{u}{v} \quad \forall \bs{x} = (y,u), \bs{\varsigma}=(w,v) \in X\,,
\end{equation}
\begin{equation}
\B{q}{\bs{x}}= a(y, q) -\scal[\Gamma_C]{u}{q} \quad \forall \bs{x}\in X, \, \forall q\in Y\,,
\end{equation}
and the functional $\mathcal{D}\in X^*$
\begin{equation}
\dual[ X^*,X]{\mathcal{D}}{\bs{\varsigma}} = \scal[{\OmegaObs}]{\dState}{w}+ \scal[\Gamma_C]{0}{v}\,,
\end{equation}
The weak form of Problem~\ref{PB:OS} has the following saddle point structure: 
\begin{problem}\label{PB:weakSP}
Find $(\bs{x},p)\in X\times Y$ such that
\begin{equation}
\begin{cases}
\A{\bs{x}}{\bs{\varsigma}} + \B{p}{\bs{\varsigma}} = \dual[ X^*,X]{\mathcal{D}}{\bs{\varsigma}} &\forall \bs{\varsigma}\in X
\\
\B{q}{\bs{x}} = \scal[\Omega]{f}{q} &\forall q\in Y
\end{cases}
\end{equation}
\end{problem}
By classical well-posedness results Problem~\ref{PB:weakSP} has a unique solution and, moreover, it is equivalent to Problem~\ref{P1:OCP}, see e.g.\ \cite{BrezziBoffiFortin}.

\section{{Numerical Discretization}}
\label{sec:VEM}
In this section, we introduce the Stabilization-Free Virtual Element (SFVEM) approximation of Problem \ref{PB:weakSP} following \cite{Berrone2024,BERRONE2025117839}. 
Let $\Mh$ denote a polygonal tessellation of $\Omega$ conforming to $\OmegaObs$. $E$ denotes a general polygon of $\Mh$ with diameter $h_E$ and let $h$ be the maximum diameter over $\Mh$. Let $e$ denote an edge of $\Mh$, $\Ih$ be the set of edges lying on $\Gamma_C$ and $\MhObs$ be the set of polygons $E$ in $\OmegaObs$.
We assume that $\Mh$ fulfills the following assumptions: $\exists \kappa>0$ such that
\begin{enumerate}[label=\textbf{A.\arabic*}]
\item\label{meshA_1} for all $E\in \Mh$, $E$ is star-shaped with respect to a
  ball of radius $\rho\geq \kappa h_E$;
\item\label{meshA_2} for all $E\in \Mh$, for all edges $e$ of $\Mh$, $\abs{e} \geq \kappa h_E$.
\end{enumerate}
For any given $E\in\Mh$, let $\Poly{k}{E}$ be the space of polynomials of degree
up to $k$ defined on $E$. Let
$\proj[\nabla]{k}{E}:\sobh{1}{E}\rightarrow\Poly{k}{E}$ be the
$\sobh{1}{}$-orthogonal projection defined up to a constant by the
orthogonality condition:
\begin{equation}\label{eq:PiNablaorthogonalitycondition}
  \forall \,w\in\sobh{1}{E}, \quad  \scal[E]{\nabla\left(\proj[\nabla]{k}{E} w -w\right)}{\nabla \pi_k}=0 \;\; \forall \, \pi_k\in\Poly{k}{E}.
\end{equation}
In order to define $\proj[\nabla]{k}{E}$ uniquely, we set
  $\displaystyle \int_{\partial E}\left(\proj[\nabla]{k}{E} w-w\right) \,\mathrm{d} s = 0$ if $k=1$ and $\displaystyle \int_{ E}\left(\proj[\nabla]{k}{E} w-w\right) \,\mathrm{d} s = 0$ otherwise.
  For any given $E\in\Mh$, let $\Vh[E]$ be the local Virtual Element Space:
\begin{multline*}
  \Vh[E] = \{ w_h\in \sobh{1}{E} \colon \Delta w_h\in \Poly{k}{E}, \, w_{h_{|e}}\in \Poly{k}{e}\,\forall e \subset \partial E, \, w_h\in\cont{\partial E},
  \\
  \scal[E]{w_h}{\pi_k} =\scal[E]{\proj[\nabla]{k}{E} w_h}{\pi_k} \; \forall \pi_k\in\Poly{k}{E}/ \Poly{k-2}{E}
  \}\,.
\end{multline*}
We recall that the degrees of freedom of this space are (see \cite{Beirao2015b})
\begin{itemize}
\item The values of $w_h$ at the vertices of $E$,
\item If $k\geq 2$, the values at $k-1$ internal distinct points on each edge $e\subset\partial E$,
\item If $k\geq 2$, the scaled moments $\frac{1}{\abs{E}}\int_E w_h\, m_i$, where $\{m_i\}_{i=1}^{\dim\Poly{k-2}{E}}$ is a basis of $\Poly{k-2}{E}$.
\end{itemize}
Moreover,
we denote the global discrete space as
\begin{equation*}
  \Vh \mydefeq \{w_h\in Y \colon w_{h_{|E}} \in \Vh[E], w_h = 0 \text{ on } \Gamma_D\}
  \,,
\end{equation*}
and we equipped it with the norm $\ennorm[Y]{\cdot}$, considering also the following notation for the local norm on each $E\in\Mh$
\begin{equation}
\ennorm[Y(E)]{w}^2= \norm[\lebl{E}]{\sqrt{\K}\nabla w}^2,
\end{equation}
for each $w\in Y(E)$.
The local approximation space for $U$ is simply defined as follows
\begin{equation}~\label{eq:polynomialDiscreteSpace}
\Uh(E) \eqdot \left\{v_h\in U: v_{h}\in \Poly{1}{e}\; \forall e\in\Ih\,, v_h\in\cont{\Gamma_C} \right\},
\end{equation}
and it is equipped with the $\lebl{}$-norm on $\Gamma_C$.
We now focus on the ingredients of the SFVEM approximation of the bilinear forms of Problem \ref{PB:weakSP}.

Let $\proj{k-1}{E}$ be the local $\lebl{}$-orthogonal projector into the space of polynomials of degree up to $k-1$, i.e. for each $w\in\lebl{E}$, $\displaystyle \int_E (\proj{k-1}{E}  w -w)\pi_{k-1}=0$ for each $\pi_{k-1}\in\Poly{k-1}{E}$.
Then, for each $\bs{w}=\left(w_1,w_2\right)^\intercal\in\lebldouble{E}$, we define the vector-valued projection operator component-wise, e.g. $\proj{k-1}{E}\bs{w}=\left(\proj{k-1}{E}w_1,\proj{k-1}{E}w_2\right)^{\intercal}$.

Moreover, under assumptions \ref{meshA_1} and \ref{meshA_2}, for each $E\in\Mh$ and for $\ell_E \in \mathbb N$, we define the space 
\begin{equation}
\mathcal{P}_{k,\ell_E} \mydefeq 
\bigl[P_{k-1}(E)\bigr]^2 \;\oplus\; 
\curl\Bigl(P_{k+\ell_E}(E) \setminus P_k(E)\Bigr)
=
\mathbf x\, P_{k-2}(E) \;\oplus\; \curl\, P_{k+\ell_E}(E),
\end{equation}
with $ \curl \pi = \left(\frac{\partial \pi}{\partial y},- \frac{\partial \pi}{\partial
    x}\right)$  for any $\pi\in \Poly{\ell_E+1}{E}$.
    Let us define $\proj{\mathcal P}{E}\nabla: H^1(E) \rightarrow \mathcal{P}_{k,\ell_E}$ the projection operator defined by the following orthogonality condition for any given $w$
\begin{equation}\label{eq:def-proj-stabfree}
  \scal[E]{\proj{\mathcal P}{E} \nabla w}{\pi} = \scal[E]{\nabla w}{\pi}  \quad\forall \pi \in \mathcal{P}_{k,\ell_E}.
\end{equation}
Notice that for $w_h \in Y_h(E)$, the polynomial projections $\proj[\nabla]{k}{E} w_h$, $\proj{k}{E} w_h$,  $\proj{k-1}{E}\nabla w_h$ and $\proj{\mathcal P}{E}\nabla w_h$ are computable exploiting only the degrees of freedom of $w_h$.

Let $\ahE{}{}:\sobh{1}{E}\times \sobh{1}{E}\to\mathbb{R}$ and $\mhE{}{}:\sobh{1}{E}\times \sobh{1}{E}\to\mathbb{R}$ be the local bilinear forms defined as, for all $w_h,q_h \in \sobh{1}{E}$
\begin{align}\label{eq:deftildeahE}
\ahE{w_h}{q_h} & \mydefeq \scal[E]{\K\proj{\mathcal P}{E} \nabla w_h}{\proj{\mathcal P}{E} \nabla q_h} \\ \nonumber
 & \qquad \qquad + \scal[E]{\beta \proj{k-1}{E}  \nabla w_h}{\proj{k}{E} q_h} + \scal[E]{\gamma \proj{k}{E}  w_h}{\proj{k}{E} q_h},
\end{align}
and 
\begin{equation}
  \mhE{w_h}{q_h} \mydefeq \scal[E]{\proj{k}{E}  w_h}{\proj{k}{E} q_h}.
\end{equation}
We define the global bilinear forms such that for each $w_h,q_h\in\sobh{1}{\Omega}$, 
\begin{equation*}
     \ah{w_h}{q_h} \mydefeq \sum_{E\in\Mh}\ahE{w_h}{q_h} \text{ and } \mh{w_h}{q_h} \mydefeq \sum_{E\in\MhObs}\mhE{w_h}{p_h}\,.
\end{equation*}

To assess the coercivity and the continuity of the bilinear forms in the SFVEM approximation, we recall the following assumption from \cite{BERRONE2025117839}
\begin{enumerate}[label=\textbf{A.\arabic*}, start=3]
\item\label{meshA_3} For any $E\in\Mh$, let $\ell_E$ be the smallest integer such that any polynomial $\pi_{k+\ell_E} \in \mathcal{P}_{k+\ell_E}(E)$ can be identified by a set of degrees of freedom
which contains $kN_E - 1$ distinct moments
$
\frac{1}{|\partial E|}\,(\pi_{k+\ell_E} ,\xi_i)_{\partial E},
$
for a scaled polynomial basis of
$ \Poly[0]{k-1}{\partial E}
\coloneqq
\left\{
\xi \;:\;
\xi|_e \in \mathbb{P}_{k-1}(e),\ \forall e \subset \partial E,
\ \int_{\partial E} \xi = 0
\right\}.
$
\end{enumerate}
We report the following result, easily derivable from \cite{BERRONE2025117839}, concerning the coercivity and the continuity of the bilinear form $a_h$, and from \cite{Beirao2015b}, concerning the continuity of $m_h$.
\begin{lemma}
\label{lem:1}
Under the assumptions \ref{meshA_1}, \ref{meshA_2}, 
for each $w_h\in\sobh{1}{\Omega}$
\begin{equation}
\mh{w_h}{w_h}\leq\norm[\lebl{\Omega}]{w_h}^2 \,,
\end{equation}
and there exists $\eta^{\ast} >0$ such that
\begin{equation} \label{eq:ah-continuity}
\ah{w_h}{w_h}\leq\eta^\ast\ennorm[Y]{w_h}^2\,.
\end{equation}
Assuming also
\ref{meshA_3}, for each $w_h\in \Vh$ there exists $\eta_{\ast}
>0$, such that, for $h$ sufficiently small, 
\begin{equation}
  \label{eq:ah-local-coercivity}
\ah{w_h}{w_h}\geq \eta_\ast\ennorm[Y]{w_h}^2 \,.
\end{equation}
\end{lemma}

\begin{remark}[Standard VEM approximation]
\label{rem:stab}
We here highlight the main difference between the proposed approximation and the one involving a standard VEM scheme (see \cite{Beirao2015b}). 
For each $E\in\Mh$, let $\atildehE{}{}:\Vh[E]\times\Vh[E]\to\mathbb{R}$ 
be the local bilinear form defined as
\begin{multline}
\label{eq:atildeVEM}
\atildehE{w_h}{q_h} \mydefeq \scal[E]{\K\proj{k-1}{E} \nabla w_h}{\proj{k-1}{E} \nabla q_h} \\ + \norm[\infty,E]{\K}\vemstab[E]{w_h-\proj[\nabla]{k}{E}w_h}{q_h-\proj[\nabla]{k}{E}q_h} \\
+ \scal[E]{\beta \proj{k-1}{E}  \nabla w_h}{\proj{k}{E} q_h}
 + \scal[E]{\gamma \proj{k}{E}  w_h}{\proj{k}{E} q_h}.
\end{multline}
The bilinear form $\vemstab[E]{}{}:\Vh[E]\times\Vh[E]\to\mathbb{R}$, called stabilization, is any symmetric positive definite bilinear form to be chosen to verify that $\exists s_\ast,s^\ast >0$ independent of $h_E$ such that
\begin{equation}
s_\ast\scal[E]{\nabla w_h}{\nabla w_h}\leq\vemstab[E]{w_h}{w_h} \leq s^\ast\scal[E]{\nabla w_h}{\nabla w_h},
\end{equation}
for each $w_h\in\Vh[E]$ verifying $\proj[\nabla]{k}{E}w_h=0$. 
    A common choice for $\vemstab[E]{}{}$ is the \textit{dofi-dofi} stabilizing bilinear form $\vemstab[E]{}{}\colon \Vh[E]\times \Vh[E]\to \mathbb{R}$, defined such that, denoting by $\chi^E(v_h)$ the vector of degrees of freedom of $v_h$ on $E$ and by $\sigma$ an arbitrary-chosen positive constant,
\begin{equation}
\vemstab[E]{w_h}{q_h}=\sigma \;\chi^E(w_h)\cdot\chi^E(q_h) \quad\forall w_h,q_h \in \Vh[E]\,. 
\end{equation}
We recall that other choices of $S^E$ can be made, allowing a relaxation of the mesh assumption \ref{meshA_2} \cite{Brenner2018}.
For each $w_h, q_h \in Y_h$, we now define the global VEM form
$$
\tilde{a}_h(w_h, q_h) \mydefeq \sum_{E\in\Mh}\atildehE{w_h}{q_h}\, ,
$$ that satisfies under the assumptions \ref{meshA_1} and \ref{meshA_2}, there exists $\tilde \eta_{\ast} \,, \tilde \eta^{\ast}  
>0$, independent of $h_E$, such that for each $w_h\in\Vh$
\begin{equation}
  \label{eq:ah-local-coercivity_VEM}
\tilde \eta_\ast\ennorm[Y]{w_h}^2\leq\tilde{a}_h(w_h, w_h)\leq \tilde \eta^\ast\ennorm[Y]{w_h}^2.
\end{equation}
\end{remark}

Let us now write the proposed approximation of Problem~\ref{PB:weakSP}.
Let us consider $\bs{x}_h \mydefeq (y_h,u_h)\in X_h$, where $X_h \mydefeq \Vh\times \Uh\subset X$ equipped with the norm $\norm[X]{\cdot}$, i.e.
$\norm[X]{\bs{x}_h}^2 = \ennorm[Y]{y_h}^2+\norm[\lebl{\Gamma_C}]{u_h}^2$ 
,
 and the bilinear forms
$\Ah{}{}:X_h\times X_h\to\mathbb{R}$
and
$\Bh{}{}:\Vh\times X_h\to\mathbb{R}$
given by
\begin{equation*}
\Ah{\bs{x}_h}{\bs{\varsigma}_h} = \mh{y_h}{w_h} + \alpha \scal[\Gamma_C]{u_h}{v_h} \quad \forall \bs{x}_h = (y_h,u_h), \bs{\varsigma}_h=(w_h,v_h) \in X_h\,,
\end{equation*}
\begin{equation}
  \label{eq:def-Bh}
\Bh{q_h}{\bs{x}_h}= \ah{y_h}{q_h} - \scal[\Gamma_C]{u_h}{q_h} \; \forall \bs{x}_h\in X_h, \, \forall q_h\in \Vh \,.
\end{equation}
and the functional $\mathcal{D}_h\in X_h^*$
\begin{equation*}
\dual[ X_h^*,X_h]{\mathcal{D}_h}{\bs{\varsigma}_h} = \sum_{E\in\MhObs}\scal[E]{\dState}{\proj{k}{E} w_h}+ \scal[\Gamma_C]{0}{v_h}\,.
\end{equation*}
Hence, we obtain the following discrete saddle point system.
\begin{problem}\label{PB:VEM-SP}
Find $(\bs{x}_h,p_h)\in X_h\times \Vh$ such that
\begin{equation}
\begin{cases}
\Ah{\bs{x}_h}{\bs{\varsigma}_h} + \Bh{p_h}{\bs{\varsigma}_h} = \dual[ X_h^*,X_h]{\mathcal{D}_h}{\bs{\varsigma}_h} &\forall \bs{\varsigma}_h\in X_h
\\
\Bh{q_h}{\bs{x}_h} =  \sum_{E\in\Mh}\scal[E]{f}{\proj{k}{E} q_h} &\forall q_h\in \Vh
\end{cases}
\end{equation}
\end{problem}
\begin{proposition}\label{prop:wellPos}
The bilinear forms
$\Ah{}{}$ and $\Bh{}{}$ are continuous and satisfy the following conditions
\begin{itemize}
\item[$\circ$] $\Ah{}{}$ is symmetric,
\item[$\circ$]$\Ah{}{}$ is coercive on $X_{h,0} \mydefeq \{\bs{\varsigma}\in X_h: \Bh{q_h}{\bs{\varsigma}_h}=0 \; \forall q_h\in \Vh\}$, i.e.\ $\exists \alpha_\ast>0$ independent of $h$ such that
\begin{equation} \label{eq:A-coercivity-on-the-kernel}
\Ah{\bs{\varsigma}_h}{\bs{\varsigma}_h} \geq \alpha_\ast\,\norm[X]{\bs{\varsigma}_h}^2 \quad \forall \bs{\varsigma}_h \in X_{h,0}.
\end{equation}
\item[$\circ$] $\Bh{}{}$ satisfies the inf-sup condition, i.e. $\exists \beta_\ast>0$ independent of $h$ such that
\begin{equation} \label{eq:B-inf-sup-cond}
\inf_{p_h\in \Vh} \sup_{\bs{\varsigma}_h\in X_h} \frac{\Bh{p_h}{\bs{\varsigma}_h}}{\ennorm[Y]{p_h}\norm[X]{\bs{\varsigma}_h}} \geq \beta_\ast\,.
\end{equation}
\end{itemize}
Then, Problem~\ref{PB:VEM-SP} admits a unique solution and is equivalent to the minimization of 
\begin{equation}
J_h(y_h,u_h) = \frac12\sum_{E\in\MhObs}\norm[\lebl{E}]{\proj{k}{E}y_h-\dState}^2+\frac{\alpha}{2}\norm[\lebl{\Gamma_C}]{u_h}^2,
\end{equation}
subject to $\ah{y_h}{q_h} = \scal[\Gamma_C]{u_h}{q_h}+\sum_{E\in\Mh}\scal[E]{f}{\proj{k}{E} q_h}$.  
\end{proposition}
\begin{proof}
It is clear that $\Ah{}{}$ is symmetric and semi-definite positive by the definitions of scalar products. Moreover, $\Ah{}{}$ is continuous, since
$\bs{x}_h = (y_h, u_h)$, $\bs{\varsigma}_h = (q_h, v_h) \in X_h$ 
\begin{align*}
|\Ah{\bs{x}_h}{\bs{\varsigma}_h}| & 
\leq \sum_{E\in\MhObs}  \norm[\lebl{E}]{\proj{k}{E}y_h}\norm[\lebl{E}]{\proj{k}{E}q_h} +  \alpha \norm[\lebl{\Gamma_C}]{u_h}\norm[\lebl{\Gamma_C}]{v_h} \\
& \leq \norm[\lebl{\Omega}]{y_h}\norm[\lebl{\Omega}]{q_h}  +  \alpha \norm[\lebl{\Gamma_C}]{u_h}\norm[\lebl{\Gamma_C}]{v_h} \\
& \leq \left ( 1 + \alpha \right )  \left(\norm[\lebl{\Omega}]{y_h}+\norm[\lebl{\Gamma_C}]{u_h}\right)\left(\norm[\lebl{\Omega}]{q_h}+\norm[\lebl{\Gamma_C}]{v_h}\right)
\\
& \leq {C^*}\left ( 1 + \alpha \right )  \norm[X]{\bs{x}_h}\norm[X]{\bs{\varsigma}_h}\,, 
\end{align*}
where we exploit Remark \ref{rem:H1norms}.
We now prove the coercivity of $\Ah{}{}$ on $X_{h,0}$. Given $\bs{\varsigma}_h \in X_{h,0}$ we have
$$
\Bh{q_h}{\bs{\varsigma}_h} = 0 \quad \forall q_h \in \Vh \quad \text{if and only if} \quad \ah{w_h}{q_h} = (v_h, q_h)_{\Gamma_C} \quad \forall q_h \in \Vh.
$$
Using Lemma \ref{lem:1}, applying the discrete trace inequality (see e.g. \cite{DiPietro2012}) and Remark \ref{rem:H1norms}, $\norm[\lebl{\Gamma_C}]{w_h}\leq C_T\ennorm[Y]{w_h}$ 
we obtain 
\begin{equation}
\begin{aligned}
\eta_\ast\ennorm[Y]{w_h}^2
\leq \ah{w_h}{w_h}
=\scal[\Gamma_C]{v_h}{w_h}
& \leq\norm[\lebl{\Gamma_C}]{v_h}\norm[\lebl{\Gamma_C}]{w_h} \\
& \leq C_T\norm[\lebl{\Gamma_C}]{v_h}\ennorm[Y]{w_h}\,,
\end{aligned}
\end{equation}
hence
\begin{equation}
\ennorm[Y]{w_h}\leq \frac{C_T}{\eta_\ast}\norm[\lebl{\Gamma_C}]{v_h}.
\end{equation}
Thus,
\begin{align*}
\Ah{\bs{\varsigma}_h}{\bs{\varsigma}_h } 
& = \sum_{E\in\MhObs}\norm[\lebl{E}]{\proj{k}{E}w_h}^2 + \alpha \norm[\lebl{\Gamma_C}]{v_h}^2 
\\
&= 
\sum_{E\in\MhObs}\norm[\lebl{E}]{\proj{k}{E}w_h}^2 + \frac{\alpha}{2} \norm[\lebl{\Gamma_C}]{v_h}^2 + \frac{\alpha}{2} \norm[\lebl{\Gamma_C}]{v_h}^2\\
								& \geq  \frac{\alpha}{2} \norm[\lebl{\Gamma_C}]{v_h}^2 + \frac{\alpha \eta_{\ast}^2}{2C^2_T} \ennorm[Y]{w_h}^2\\
&  \geq \frac{\alpha}{2}\min \left ( 1, \frac{\eta_{\ast}^2}{C^2_T}\right )  \left( \ennorm[Y]{w_h}^2 +  \norm[\lebl{\Gamma_C}]{v_h}^2 \right ) 
\geq  \frac{\alpha}{2}\min \left ( 1, \frac{\eta_{\ast}^2}{C^2_T}\right ) \norm[X]{\bs{\varsigma}_h}^2.
\end{align*}
Furthermore, $\Bh{}{}$ is continuous on $\Vh \times X_h$ by definition. Moreover, $\Bh{}{}$ is inf-sup stable, indeed applying Lemma \ref{lem:1} we have:
\begin{align*}
\sup_{0 \neq \bs{x}_h \in X_h} \frac{\Bh{q_h}{\bs{x}_h }}{\norm[X]{\bs{x}_h}} 
& = \sup_{0 \neq  (y_h, u_h) \in \Vh \times U_h} \frac{\ah{y_h}{q_h} - (u_h, q_h)_{\Gamma_C}}{\sqrt{\ennorm[Y]{{y}_h}^2 + 
\norm[U_h]{{v}_h}^2}}\\
& \geq  \frac{\ah{q_h}{q_h}}{\ennorm[Y]{{q}_h}} \geq \eta_{\ast} \ennorm[Y]{{q}_h},
\\
\end{align*}
choosing $(y_h, u_h)$ as $(q_h, 0)$ to bound the supremum from below. The arbitrariness of $q_h$ concludes the proof.
\end{proof}

\section{Error Estimates}
\label{sec:err_estimates}
In order to derive optimal a priori error estimates, we have to provide the polynomial approximation properties and the interpolation error estimates, presented in Lemmas \ref{lem:polynomial_approximation}, \ref{lem:higher_order_interpolationError} and \ref{lem:Uh-interpolant}, respectively. 
The results about polynomial approximation can be led back to classical results,
while the VEM interpolation error estimate is proved in \cite{MRR2015} and the interpolation on $\Uh$ derives from standard results of one-dimensional FEM interpolation.

\begin{lemma} \label{lem:polynomial_approximation}
Under assumption \ref{meshA_1}, given  
$\proj{k}{E}$ be the local $\lebl{}$-orthogonal projection into $\Poly{k}{E}$ and $\proj[\nabla]{k}{E}$ $\sobh{1}{}$-orthogonal projection defined by \eqref{eq:PiNablaorthogonalitycondition},
then, there exists $ C>0$, independent of $h_E$, such that $\forall \phi\in\sobh{s}{E}$
\begin{equation} \label{eq:L2_polynomial_approximation}
\seminorm[m,E]{\phi-\proj{k}{E} \phi} \leq C h_{E}^{s-m} \seminorm[s,E]{\phi}, \; \mbox{ with } m\leq s \leq k+1 \, ,
\end{equation} 
and if $s\geq 1$
\begin{equation}\label{eq:H1_polynomial_approximation}
\seminorm[m,E]{\phi-\proj[\nabla]{k}{E} \phi} \leq C h_{E}^{s-m} \seminorm[s,E]{\phi}, \; \mbox{ with } m\leq s \leq k+1 \,.
\end{equation} 
\end{lemma}

\begin{lemma} \label{lem:higher_order_interpolationError}
Under assumptions \ref{meshA_1} and \ref{meshA_2}, let 
$w\in\sobh{s+1}{\Omega}$,
$0\leq s\leq k$, 
there exists $C>0$, independent of $h$, such that $\forall \,h$, $\exists w_I\in\Vh$ satisfying
\begin{equation}\label{eq:higher_order_interpolationError}
\norm[\lebl{\Omega}]{w-w_I}+h\norm[\lebl{\Omega}]{\nabla(w-w_I)} \leq C h^{s+1}\seminorm[s+1,\Omega]{w}.
\end{equation}
Thus,
\begin{equation}\label{eq:Ynorm_interpEstimate}
\ennorm[Y]{w-w_I} \leq C h^{s}\seminorm[s+1,\Omega]{w}.
\end{equation}
\end{lemma}

\begin{lemma}
  \label{lem:Uh-interpolant}
  Under assumptions \ref{meshA_1} and \ref{meshA_2}, let
  
  $v\in \sobh{r}{\Gamma_C}$, 
  $0\leq r \leq k+1$. 
  Then there exists a function
  $v_I\in\Uh$ and a constant $C>0$ independent of $h$ such that
  \begin{equation}
    \label{eq:Uh-interpolant}
    \norm[\lebl{\Gamma_C}]{v-v_I} \leq C h^r \seminorm[r,\Gamma_C]{v} \,.
  \end{equation}
\end{lemma}
In order to prove Theorem \ref{thm:ErrorEstimate}, we need this additional auxiliary result.
\begin{lemma}
  \label{lem:problem-interpolant}
  Let $((y,u),p)\in X\times Y$ be the solution to Problem \ref{PB:weakSP} and
  let $((y_h,u_h),p_h)\in X_h\times \Vh$ be the solution to Problem
  \ref{PB:VEM-SP}. Moreover, let $u_I\in U_h$ as in Lemma
  \ref{lem:Uh-interpolant}. Then there esists $y_{\Bh{}{}} \in \Vh$ such
  that
  \begin{equation}
    \label{eq:def-yBh}
    \Bh{q_h}{\left(y_{\Bh{}{}}-y_h,u_I - u_h\right)} = 0
    \quad \forall q_h\in\Vh \,.
  \end{equation}
  Moreover, if $y\in\sobh{s+1}{\Omega}$, $u\in\sobh{r}{\Gamma_C}$,
  $f\in\sobh{t}{\Omega}$ with $0\leq s \leq k$, $0\leq r \leq k+1$ and
  $0\leq t \leq k+1$, then there exists $C>0$ independent of $h$ such that
  \begin{equation}
    \label{eq:estim-yBh}
    \ennorm[Y]{y - y_{\Bh{}{}}} \leq
    C \left(h^s \seminorm[s+1,\Omega]{y} + h^{t+1}\seminorm[t,\Omega]{f} +
    h^{r} \seminorm[r,\Gamma_C]{u}\right) \,.
  \end{equation}
\end{lemma}

\begin{proof}
  Let $y_{\Bh{}{}}\in\Vh$ be the unique solution to
  \begin{equation}
    \label{eq:prob-yBh}
    \ah{y_{\Bh{}{}}}{q_h} = \scal[\Gamma_C]{u_I}{q_h}
    + \sum_{E\in\Mh} \scal{f}{\proj{k}{E} q_h} \quad \forall q_h \in \Vh \,.
  \end{equation}
  Notice that, by definition of $\Bh{}{}$ \eqref{eq:def-Bh} and Problem \ref{PB:VEM-SP}, we have
  \begin{equation*}
    \Bh{q_h}{(y_{\Bh{}{}},u_I)} = \sum_{E\in\Mh} \scal{ f}{\proj{k}{E} q_h} =
    \Bh{q_h}{(y_h,u_h)}\,,
  \end{equation*}
  thus \eqref{eq:def-yBh} holds true. 
  Moreover, let $y_I\in\Vh$ be the function
  satisfying \eqref{eq:Ynorm_interpEstimate}, then by a triangle
  inequality we have
  \begin{equation*}
    \ennorm[Y]{y - y_{\Bh{}{}}}
    \leq \ennorm[Y]{y - y_I} + \ennorm[Y]{y_I - y_{\Bh{}{}}}
    \leq C h^s \seminorm[s+1,\Omega]{y} +
    \ennorm[Y]{y_I - y_{\Bh{}{}}} \,.
  \end{equation*}
  Regarding the second term above, let
  $e_h = y_I-y_{\Bh{}{}}$. Exploiting \eqref{eq:ah-local-coercivity},
  \eqref{eq:prob-yBh} and the second equation of Problem \ref{PB:weakSP} we get
  \begin{equation}
    \label{eq:yI-yBh-estim}
    \begin{split}
      \ennorm[Y]{y_I - y_{\Bh{}{}}}^2
      &\leq C \ah{y_I - y_{\Bh{}{}}}{e_h}
      \\
      &=
        C\Big (\ah{y_I}{e_h} - \a{y}{e_h} \\
        & \qquad \qquad + \sum_{E\in\Mh} \scal[E]{f-\proj{k}{E}f}{e_h}
        + \scal[\Gamma_C]{u-u_I}{e_h}\Big ) \,.
    \end{split}
  \end{equation}
  The terms above are estimated separately, as follows. 
  Adding and subtracting $\sum_{E\in\Mh}\ahE{y}{e_h}$, applying the continuity of $\ah{}{}$ \eqref{eq:ah-continuity} and the estimate \eqref{eq:Ynorm_interpEstimate}, we get
  \begin{equation}\label{eq:estim-ahInterp-aSol}
      \begin{split}
          \ah{y_I}{e_h} - \a{y}{e_h}
      &= \sum_{E\in\Mh} \ahE{y_I - y}{e_h}
        + \ahE{y}{e_h} - \aE{y}{e_h} \\
        & \leq Ch^s \seminorm[s+1,\Omega]{y} \ennorm[Y]{e_h} +
      \sum_{E\in\Mh}\ahE{y}{e_h} - \aE{y}{e_h}
      \\
      & \leq Ch^s \seminorm[s+1,\Omega]{y} \ennorm[Y]{e_h}\,,
      \end{split}
  \end{equation}
  where the last estimate can be derived as in \cite[Theorem 5]{BERRONE2025117839}.
  
  Next, the term of \eqref{eq:yI-yBh-estim} depending on $f$ is estimated
  exploiting the definition of $\proj{k}{E}$, a Cauchy-Schwarz
  inequality and estimate \eqref{eq:L2_polynomial_approximation}:
  \begin{equation*}
    \begin{split}
      \sum_{E\in\Mh}\scal[E]{f-\proj{k}{E}f}{e_h}
      &=
        \sum_{E\in\Mh}\scal[E]{f-\proj{k}{E}f}{e_h - \proj{k}{E}e_h}
      \\
      &\leq
        \sum_{E\in\Mh}\norm[\lebl{E}]{f-\proj{k}{E}f} \norm[\lebl{E}]{e_h - \proj{k}{E}e_h}
      \\
      &\leq C h^{t} \seminorm[t,\Omega]{f} \cdot h \norm[\lebl{\Omega}]{\nabla e_h}
        \leq C h^{t+1} \seminorm[t,\Omega]{f} \ennorm[Y]{e_h} \,.
    \end{split}
  \end{equation*}
  Finally, the last term in \eqref{eq:yI-yBh-estim} is bounded by a
  Cauchy-Schwarz inequality, the estimate \eqref{eq:Uh-interpolant} and a trace
  inequality:
  \begin{equation*}
    \begin{split}
      \scal[\Gamma_C]{u-u_I}{e_h}
      &\leq \norm[\lebl{\Gamma_C}]{u-u_I} \norm[\lebl{\Gamma_C}]{e_h}
        \leq C h^r \seminorm[r,\Gamma_C]{u} \norm[\lebl{\Gamma_C}]{e_h} \\ 
        & \leq
        C h^{r} \seminorm[r,\Gamma_C]{u} \norm[\lebl{\Omega}]{\nabla e_h}
      \\
      &\leq C h^{r} \seminorm[r,\Gamma_C]{u} \ennorm[Y]{e_h} \,.
    \end{split}
  \end{equation*}
  Using the above estimates in \eqref{eq:yI-yBh-estim} and simplifying
  $\ennorm[Y]{e_h}$, we get \eqref{eq:estim-yBh}.
\end{proof}
\begin{theorem} \label{thm:ErrorEstimate}
Let 
$(\bs{x},p)=((y,u),p)\in \left(\sobh{s+1}{\Omega}\cap\sobh[0]{1}{\Omega}\right) \times  \sobh{r}{\Gamma_C} \times \left(\sobh{s+1}{\Omega}\cap\sobh[0]{1}{\Omega}\right)$ be the solution to Problem \ref{PB:weakSP}, with $0\leq s\leq k$ and $0\leq r \leq k+1$. 
Let $f\in \sobh{t}{\Omega}$ and $\dState \in \sobh{\tau}{\OmegaObs}$ be the right-hand side of Problem \ref{PB:weakSP}, with $0\leq t\leq k+1$ and $0\leq\tau\leq k+1$.
Then, there exists a constant $C$ independent of $h$ such that the unique solution $((y_h,u_h),p_h)\in X_h\times \Vh$ to Problem
  \ref{PB:VEM-SP}
 satisfies the following error estimate:
 \begin{equation}\label{eq:ErrorEstimate}
 \begin{split}
\norm[X]{\bs{x}-\bs{x}_h}+\ennorm[Y]{p-p_h}
 & \leq
    C \big (
    h^s \seminorm[s+1,\Omega]{y}  +
    h^{r} \seminorm[r,\Gamma_C]{u} 
    + h^s \seminorm[s+1,\Omega]{p}\\
   & \qquad \qquad \qquad  + h^{t+1}\seminorm[t,\Omega]{f}
    + h^{\tau+1}\seminorm[\tau,\OmegaObs]{\dState}
    \big )
    \end{split}
\end{equation}
\end{theorem}
\begin{proof}
First, let us prove the estimate for $\norm[X]{\bs{x}-\bs{x}_h}$.
Let $\bs{x}_I\mydefeq(y_{\Bh{}{}},u_I)$ be given by Lemma \ref{lem:Uh-interpolant} and \ref{lem:problem-interpolant}.
Applying the triangle inequality, we have
\begin{equation}\label{eq:triang-XI}
\norm[X]{\bs{x}-\bs{x}_h} \leq \norm[X]{\bs{x}-\bs{x}_I} +\norm[X]{\bs{x}_I-\bs{x}_h}\,,
\end{equation}
then the term $\norm[X]{\bs{x}-\bs{x}_I}$ can be estimated applying \eqref{eq:estim-yBh} and \eqref{eq:Uh-interpolant}, i.e.
\begin{equation}\label{eq:estim-X-XI}
\begin{split}
    \norm[X]{\bs{x}-\bs{x}_I} & = 
    \ennorm[Y]{y - y_{\Bh{}{}}} + \norm[\lebl{\Gamma_C}]{u-u_I}\\
   & \lesssim
     \left(h^s \seminorm[s+1,\Omega]{y} + h^{t+1}\seminorm[t,\Omega]{f} +
    h^{r} \seminorm[r,\Gamma_C]{u}\right)
\end{split}
\end{equation}
where $ \xi \lesssim \zeta$ denotes that there exists a constant $C$ independent of $h$ such that $ \xi \leq C \zeta$. This notation is widely used in the following.
On the other hand, let $\bs{\varepsilon}_h=\bs{x}_I-\bs{x}_h$.  Noticing that $\bs{\varepsilon}_h\in X_{h,0}$ thanks to Lemma \ref{lem:problem-interpolant}, we can apply the coercivity of $\Ah{}{}$ on $X_{h,0}$  \eqref{eq:A-coercivity-on-the-kernel} and, adding and subtracting terms, we obtain
\begin{equation}\label{eq:1-est-X}
\begin{aligned}
\alpha_\ast \norm[X]{\bs{\varepsilon}_h}^2 
& \leq \Ah{\bs{x}_I-\bs{x}_h}{\bs{\varepsilon}_h} \\
 & = \Ah{\bs{x}_I-\bs{x}}{\bs{\varepsilon}_h} +\Ah{\bs{x}}{\bs{\varepsilon}_h} \mp \A{\bs{x}}{\bs{\varepsilon}_h} -\Ah{\bs{x}_h}{\bs{\varepsilon}_h}\,.
\end{aligned}
\end{equation}
Exploiting the continuity of $\Ah{}{}$ and \eqref{eq:estim-X-XI}, the first term of \eqref{eq:1-est-X} can be estimated as
\begin{equation}\label{eq:1-term-est-X}
\Ah{\bs{x}_I-\bs{x}}{\bs{\varepsilon}_h} \lesssim \left(h^s \seminorm[s+1,\Omega]{y} + h^{t+1}\seminorm[t,\Omega]{f} +
    h^{r} \seminorm[r,\Gamma_C]{u}\right) \norm[X]{\bs{\varepsilon}_h}\,.
\end{equation}
Then, the term $\Ah{\bs{x}}{\bs{\varepsilon}_h} - \A{\bs{x}}{\bs{\varepsilon}_h}$ of \eqref{eq:1-est-X} can be bounded exploiting the definition of $\proj{k}{E}$, the Cauchy-Schwarz inequality and \eqref{eq:L2_polynomial_approximation}, i.e.
\begin{equation}\label{eq:2-term-est-X}
\begin{aligned}
\Ah{\bs{x}}{\bs{\varepsilon}_h} - \A{\bs{x}}{\bs{\varepsilon}_h}
&= \mh{y}{y_{\Bh{}{}}-y_h} - \scal[\OmegaObs]{y}{y_{\Bh{}{}}-y_h} 
\\
&=\sum_{E\in\MhObs} \scal[E]{\proj{k}{E}  y}{\proj{k}{E}(y_{\Bh{}{}}-y_h)} \\
& \qquad \qquad \qquad \qquad  - \sum_{E\in\MhObs} \scal[E]{y}{(y_{\Bh{}{}}-y_h)}
\\
&= \sum_{E\in\MhObs} \scal[E]{\proj{k}{E}  y-y}{(y_{\Bh{}{}}-y_h) - \proj{k}{E}(y_{\Bh{}{}}-y_h)}
\\
&\lesssim h^{s+1}\seminorm[s+1,\Omega]{y} h \norm[\lebl{\Omega}]{\nabla (y_{\Bh{}{}}-y_h)} \\
& \lesssim h^{s+2}\seminorm[s+1,\Omega]{y} \norm[X]{\bs{\varepsilon}_h}\,.
\end{aligned}
\end{equation}
The remaining term in \eqref{eq:1-est-X}, $\A{\bs{x}}{\bs{\varepsilon}_h} -\Ah{\bs{x}_h}{\bs{\varepsilon}_h}$, can be rewritten as
\begin{equation}\label{eq:3-term-est-X}
\begin{aligned}
\A{\bs{x}}{\bs{\varepsilon}_h} -\Ah{\bs{x}_h}{\bs{\varepsilon}_h}
&= \Bh{p_h}{\bs{\varepsilon}_h}-\B{p}{\bs{\varepsilon}_h} \\
& \quad \quad + 
\sum_{E\in\MhObs} \scal[E]{  \dState}{(y_{\Bh{}{}}-y_h)-\proj{k}{E}(y_{\Bh{}{}}-y_h)}
\\
&= \Bh{p_I}{\bs{\varepsilon}_h}-\B{p}{\bs{\varepsilon}_h} \\
& \qquad \qquad +
\sum_{E\in\MhObs} \scal[E]{  \dState-\proj{k}{E} \dState}{(y_{\Bh{}{}}-y_h)}
\\
&= \ah{y_{\Bh{}{}}-y_h}{p_I}-\a{y_{\Bh{}{}}-y_h}{p}\\ 
& \qquad \qquad +\scal[\Gamma_C]{u_I-u_h}{p-p_I}
\\
& \qquad \qquad \qquad
+
\sum_{E\in\MhObs} \scal[E]{  \dState-\proj{k}{E}\dState}{(y_{\Bh{}{}}-y_h)}
\end{aligned}
\end{equation}
where we apply Problems \ref{PB:weakSP} and \ref{PB:VEM-SP}, Lemma \ref{lem:higher_order_interpolationError} to define $p_I\in\Vh$ interpolating $p$, Lemma \ref{lem:problem-interpolant} and the definition of $\proj{k}{E}$.
Then, the term $\ah{y_{\Bh{}{}}-y_h}{p_I}-\a{y_{\Bh{}{}}-y_h}{p}$ is estimated as
\begin{equation}
    \ah{y_{\Bh{}{}}-y_h}{p_I}-\a{y_{\Bh{}{}}-y_h}{p} \lesssim h^s \seminorm[s+1,\Omega]{p}\ennorm[Y]{y_{\Bh{}{}}-y_h}\,,
\end{equation}
following the same steps applied in the Proof of Lemma \ref{lem:problem-interpolant} (see \eqref{eq:estim-ahInterp-aSol}).
The term $\scal[\Gamma_C]{u_I-u_h}{p-p_I}$ is estimated using the Cauchy-Schwarz inequality, a discrete trace inequality together with a Poincaré inequality and \eqref{eq:higher_order_interpolationError}, obtaining
\begin{equation}
    \scal[\Gamma_C]{u_I-u_h}{p-p_I} \lesssim h^{s}\seminorm[s+1,\Omega]{p}\norm[\lebl{\Gamma_C}]{u_I-u_h}\,. 
\end{equation}
The last term of \eqref{eq:3-term-est-X} is estimated using the Cauchy-Schwarz inequality  and \eqref{eq:L2_polynomial_approximation} as
\begin{equation}
\begin{aligned}
    & \sum_{E\in\MhObs} \scal[E]{  \dState-\proj{k}{E}\dState}{(y_{\Bh{}{}}-y_h)} \\
    & \qquad \qquad = \sum_{E\in\MhObs} \scal[E]{ \dState-\proj{k}{E}\dState}{(y_{\Bh{}{}}-y_h)-\proj{k}{E}(y_{\Bh{}{}}-y_h)} 
    \\
    & \qquad \qquad \lesssim h^\tau \seminorm[\tau,\OmegaObs]{\dState} h \norm[\lebl{\Omega}]{\nabla(y_{\Bh{}{}}-y_h)}
    \\
    &\qquad \qquad  \lesssim h^{\tau+1} \seminorm[\tau,\OmegaObs]{\dState}\ennorm[Y]{y_{\Bh{}{}}-y_h}\,.
    \end{aligned}
\end{equation}
Hence, coming back to \eqref{eq:3-term-est-X}, we obtain that
\begin{equation}\label{eq:3-term-est-X-final}
    \A{\bs{x}}{\bs{\varepsilon}_h} -\Ah{\bs{x}_h}{\bs{\varepsilon}_h} \lesssim \left(h^s \seminorm[s+1,\Omega]{p} + h^{\tau+1} \seminorm[\tau,\OmegaObs]{\dState} \right) \norm[X]{\bs{\varepsilon}_h}\,.
\end{equation}
Substituting in \eqref{eq:1-est-X} the estimates \eqref{eq:1-term-est-X}, \eqref{eq:2-term-est-X} and \eqref{eq:3-term-est-X-final}, we get
\begin{equation}
    \begin{split}    
    \norm[X]{\bs{x}_I-\bs{x}_h} & \lesssim 
    h^s \seminorm[s+1,\Omega]{y}  +
    h^{r} \seminorm[r,\Gamma_C]{u} \\
    & \qquad \qquad \qquad 
    + h^s \seminorm[s+1,\Omega]{p}
    + h^{t+1}\seminorm[t,\Omega]{f}+ h^{\tau+1} \seminorm[\tau,\OmegaObs]{\dState} \,.
\end{split}  
\end{equation}
Hence, considering also \eqref{eq:triang-XI} and \eqref{eq:estim-X-XI}, we have the estimate for $\norm[X]{\bs{x}-\bs{x}_h}$, i.e.
\begin{equation}\label{eq:ESTIMATE-X-Xh}
\begin{split}
        \norm[X]{\bs{x}-\bs{x}_h} & \lesssim     h^s \seminorm[s+1,\Omega]{y}  +
    h^{r} \seminorm[r,\Gamma_C]{u} \\
     & \qquad \qquad \qquad
    + h^s \seminorm[s+1,\Omega]{p}
    + h^{t+1}\seminorm[t,\Omega]{f}+ h^{\tau+1} \seminorm[\tau,\OmegaObs]{\dState} \,.
\end{split}
\end{equation}
Now, let us focus on the estimate of $\ennorm[Y]{p-p_h}$.
Let $p_I\in\Vh$, interpolating $p$, be given by Lemma \ref{lem:higher_order_interpolationError}. Applying the triangle inequality and \eqref{eq:Ynorm_interpEstimate}, we have
\begin{equation}
    \ennorm[Y]{p-p_h}\leq \ennorm[Y]{p-p_I} +\ennorm[Y]{p_I-p_h}\lesssim h^s \seminorm[s+1,\Omega]{p}  +\ennorm[Y]{p_I-p_h}\,.
\end{equation}
On the other hand, starting from the inf-sup condition \eqref{eq:B-inf-sup-cond}, we have
\begin{equation}\label{eq:inf-sup-step-est-P}
    \ennorm[Y]{p_I-p_h} \leq \frac{1}{\beta_\ast} \sup_{\bs{\varsigma}_h\in X_h} \frac{\Bh{p_I-p_h}{\bs{\varsigma}_h}}{\norm[X]{\bs{\varsigma}_h}}\,.
\end{equation}
Let us estimate $\Bh{p_I-p_h}{\bs{\varsigma}_h}$ for any given $\bs{\varsigma}_h\in X_h$.
First, adding and subtracting $\B{p}{\bs{\varsigma}_h}$, defining $\bs{\varsigma}_h=(w_h,v_h)$ we get
\begin{equation}\label{eq:1-est-P}
    \begin{aligned}
        \Bh{p_I-p_h}{\bs{\varsigma}_h} &= 
        \Bh{p_I}{\bs{\varsigma}_h} \mp \B{p}{\bs{\varsigma}_h}
        -\Bh{p_h}{\bs{\varsigma}_h} 
        \\
        &= \ah{w_h}{p_I}-\a{w_h}{p}  + \scal[\Gamma_C]{v_h}{ p-p_I} \\
        & \qquad \qquad \qquad \qquad \qquad
        +\B{p}{\bs{\varsigma}_h} -\Bh{p_h}{\bs{\varsigma}_h}
        \,.
    \end{aligned}
\end{equation}
{The first two terms are estimated similarly as in Lemma \ref{lem:problem-interpolant} (see \eqref{eq:estim-ahInterp-aSol}), hence adding and subtracting $\sum_{E\in\Mh}\ahE{w_h}{p}$, applying the continuity of $\ah{}{}$ \eqref{eq:ah-continuity} and the estimate \eqref{eq:Ynorm_interpEstimate}, we get 
\begin{equation} \label{eq:1-term-est-P}
\begin{aligned}    
\ah{w_h}{p_I}-\a{w_h}{p}  
&= \sum_{E\in\Mh}\ahE{w_h}{p_I-p} +\ahE{w_h}{p}-\aE{w_h}{p} 
\\
&\lesssim  h^s \seminorm[s+1,\Omega]{p}\ennorm[Y]{w_h} + \sum_{E\in\Mh}\ahE{w_h}{p}-\aE{w_h}{p}
\\
&\lesssim  h^s \seminorm[s+1,\Omega]{p}\ennorm[Y]{w_h}\,,
\end{aligned}
\end{equation}
where the last estimate can be derived for the diffusion and the reaction term, thanks to their symmetry, as in \cite[Theorem 5]{BERRONE2025117839}. For what concerns the transport term, we derive the estimate as in \cite[Theorem 3]{MMConvDiff}.
}
Then, applying the Cauchy-Schwarz inequality, a discrete trace inequality together with a Poincaré inequality and \eqref{eq:higher_order_interpolationError},  the term $\scal[\Gamma_C]{v_h}{ p-p_I}$ in \eqref{eq:1-est-P} is bounded as
\begin{equation}\label{eq:2-term-est-P}
\begin{split}
        \scal[\Gamma_C]{v_h}{ p-p_I} & \leq \norm[\lebl{\Gamma_C}]{v_h}\norm[\lebl{\Gamma_C}]{p-p_I}
         \lesssim h^{s} \seminorm[s+1,\Omega]{p} \norm[\lebl{\Gamma_C}]{v_h}\,.
\end{split}
\end{equation}
The last two terms of \eqref{eq:1-est-P}, applying Problems \ref{PB:weakSP} and \ref{PB:VEM-SP}, become
\begin{equation}\label{eq:3-term-est-P}
\begin{split}
        \B{p}{\bs{\varsigma}_h}-\Bh{p_h}{\bs{\varsigma}_h} & = -\A{\bs{x}}{\bs{\varsigma}_h} +\scal[{\OmegaObs}]{\dState}{w_h} + \Ah{\bs{x}_h}{\bs{\varsigma}_h} \\ 
        & \qquad \qquad \qquad \qquad
        - \sum_{E\in\MhObs}\scal[E]{\dState}{\proj{k}{E} w_h}\,.
\end{split}
\end{equation}
Applying the Cauchy-Schwarz inequality and \eqref{eq:L2_polynomial_approximation}, we get
\begin{multline}\label{eq:3p1-term-est-P}
    \sum_{E\in\MhObs}\scal[E]{\dState}{w_h-\proj{k}{E} w_h}  = \sum_{E\in\MhObs}\scal[E]{\dState-\proj{k}{E}\dState}{w_h-\proj{k}{E} w_h}
    \\
     \leq \sum_{E\in\MhObs}\norm[\lebl{E}]{\dState-\proj{k}{E}\dState}\norm[\lebl{E}]{w_h-\proj{k}{E} w_h}
    \\
     \lesssim h^\tau \seminorm[\tau,\OmegaObs]{\dState} h \norm[\lebl{\Omega}]{\nabla w_h} \lesssim h^{\tau+1} \seminorm[\tau,\OmegaObs]{\dState} \ennorm[Y]{w_h}
    \,.
\end{multline}
On the other hand, the term $\Ah{\bs{x}_h}{\bs{\varsigma}_h} -\A{\bs{x}}{\bs{\varsigma}_h}$ in \eqref{eq:3-term-est-P} is estimated adding and subtracting $\Ah{\bs{x}}{\bs{\varsigma}_h}$, applying the definition of $\Ah{}{}$ and $\A{}{}$, the definition and the continuity of $\proj{k}{E}$ together with the Cauchy-Schwarz inequality, a Poincaré inequality and \eqref{eq:L2_polynomial_approximation}, and \eqref{eq:ESTIMATE-X-Xh}, i.e.
\begin{equation}\label{eq:3p2-term-est-P}
    \begin{aligned}
        \Ah{\bs{x}_h}{\bs{\varsigma}_h} -\A{\bs{x}}{\bs{\varsigma}_h} 
        &= \Ah{\bs{x}_h-\bs{x}}{\bs{\varsigma}_h} +\Ah{\bs{x}}{\bs{\varsigma}_h} -\A{\bs{x}}{\bs{\varsigma}_h}
        \\
        &= \sum_{E\in\MhObs} \scal[E]{\proj{k}{E}  (y_h-y)}{\proj{k}{E} w_h} \\
        & \quad + \alpha \scal[\Gamma_C]{u_h-u}{v_h}
        \\
        & \quad \quad +\sum_{E\in\MhObs} \scal[E]{\proj{k}{E}  y}{\proj{k}{E} w_h} - \scal[E]{y}{w_h}
        \\
        &\leq \norm[\lebl{\Omega}]{y-y_h}\norm[\lebl{E}]{w_h} \\
        & \quad + \alpha \norm[\lebl{\Gamma_C}]{u-u_h}\norm[\lebl{\Gamma_C}]{v_h}
        \\
        &\quad +\sum_{E\in\MhObs} \scal[E]{y-\proj{k}{E}  y}{\proj{k}{E} w_h-w_h}
        \\
        &\lesssim \norm[X]{\bs{x}-\bs{x}_h}\norm[X]{\bs{\varsigma}_h} + h^{s+1}\seminorm[s+1]{y} \ennorm[Y]{w_h}
        \\
        &\lesssim\left( h^s \seminorm[s+1,\Omega]{y}  +
    h^{r} \seminorm[r,\Gamma_C]{u} 
    + h^s \seminorm[s+1,\Omega]{p}\right)\norm[X]{\bs{\varsigma}_h}\, \\
    & \quad \quad 
    + \left(h^{t+1}\seminorm[t,\Omega]{f}+ h^{\tau+1} \seminorm[\tau,\OmegaObs]{\dState}\right)\norm[X]{\bs{\varsigma}_h}\,.
    \end{aligned}
\end{equation}
Finally, considering together \eqref{eq:1-term-est-P}, \eqref{eq:2-term-est-P}, \eqref{eq:3p1-term-est-P} and \eqref{eq:3p2-term-est-P}, and substituting in \eqref{eq:inf-sup-step-est-P}, we obtain the estimate
\begin{equation*}
    \ennorm[Y]{p_I-p_h}\lesssim h^s \seminorm[s+1,\Omega]{y}  +
    h^{r} \seminorm[r,\Gamma_C]{u} 
    + h^s \seminorm[s+1,\Omega]{p}
    + h^{t+1}\seminorm[t,\Omega]{f}+ h^{\tau+1} \seminorm[\tau,\OmegaObs]{\dState}\,,
\end{equation*}
that together with \eqref{eq:ESTIMATE-X-Xh} gives the thesis \eqref{eq:ErrorEstimate}.
\end{proof}

\begin{remark}[Standard VEM approximation]
Substituting $a_h$ with $\tilde a_h$ defined in \eqref{eq:atildeVEM}, Theorem \ref{thm:ErrorEstimate} still holds true, following \cite{Beirao2015b} for the VEM counterpart of estimates \eqref{eq:estim-ahInterp-aSol} and \eqref{eq:1-term-est-P}.
To the best of the authors' knowledge, the complete a priori error estimates of a Neumann boundary optimal control problem in a saddle point framework, have never been addressed before in the VEM community.
\end{remark}

\section{Numerical Results}
\label{sec:results}
{
In this section, we present three numerical experiments designed to assess the stability and robustness of the proposed method. The first test is a convergence study aimed at validating the theoretical results established in the previous section. The second test examines the dependence of the standard VEM procedure—described in Remark~\ref{rem:stab}—on the choice of the stabilization parameter~$\sigma$, and compares its accuracy with that of the proposed method. Finally, the third test introduces a more application‑oriented setting to assess the robustness of the proposed approach.
}

\subsection{Test 1}\label{sec:test1}
\begin{figure}
  \begin{center}
    \begin{tikzpicture}[scale=3]

      \filldraw[color=black, very thick, dashed](0,0) -- (0,1);
      \filldraw[color=black, very thick](0,0) -- (1,0); \filldraw[color=black,
      very thick](1,0.) -- (1., 1.); \filldraw[color=black, very thick]((1.,1.)
      -- (0,1);

      \node at (0,-.1){\color{black}{$(0,0)$}}; \node at
      (1,-.1){\color{black}{$(1,0)$}}; \node at (1,1.1){\color{black}{$(1,1)$}};
      \node at (0,1.1){\color{black}{$(0,1)$}}; \node at
      (.5,.5){\color{black}{$\Omega$}}; \node at
      (1.15,.5){\color{black}{$\Gamma_D$}}; \node at
      (-.15,.5){\color{black}{$\Gamma_C$}};

    \end{tikzpicture}
  \end{center}
  \caption{ Test
    1: schematic representation of the domain.  
    }
  \label{fig:domain}
\end{figure}
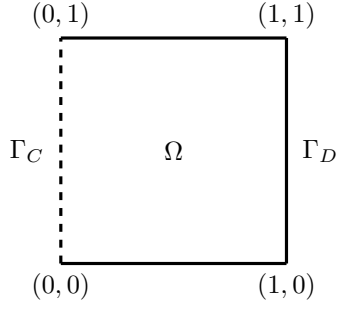

The convergence test is carried out on a $\Omega=(0,1)^2$. The boundary is
partitioned in $\Gamma_N = \emptyset$, $\Gamma_C = \{0\} \times (0,1)$ and
$\Gamma_D = \partial \Omega \setminus \Gamma_C$. The spatial coordinates are
denoted with $x = (x_1, x_2)$. Figure \ref{fig:domain} depicts a schematic representation of the domain. In this case, we deal with
\emph{distributed observation}, i.e.\ $\OmegaObs=\Omega$.

Let 
\begin{equation*}
  y(x_1,x_2) = (1 - x_1)\sin(\pi x_2), \quad u(x_1,x_2) = \sin(\pi x_2) 
  \end{equation*}
  and
  $$
  p(x_1,x_2) = \sin(\pi x_2)\cos \left (\frac{\pi}{2}x_1 \right )
$$
  be the exact solutions of Problem \ref{PB:weakSP} with 
 $\alpha =1$,
$\K(x)  \equiv \gamma(x) \equiv 1$, $\beta = [0,1]$ and
{
\begin{align*}
  f(x_1,x_2)
  &= (1-x_1)
  \left [ \sin(\pi x_2)(1 + \pi^2) + \pi \cos(\pi x_2)\right ]\,,
  \\
  \dState(x_1,x_2)
  &= \sin(\pi x_2)\left [ (1 - x_1) +  \left ( \frac{5\pi^2}{4} +  1\right ) \cos \left
    (\frac{\pi}{2}x_1 \right ) \right ] - \pi \cos(\pi x_2) \cos \left (\frac{\pi}{2}x_1\right ).
\end{align*}
}
\begin{figure}
  \centering
  \begin{subfigure}[t]{0.33\textwidth}
    \centering \includegraphics[scale=0.267]{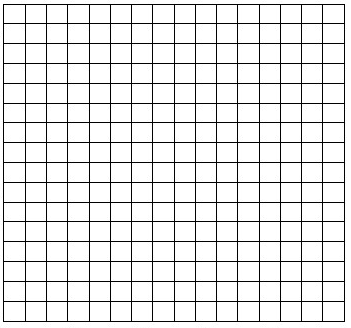}
    \caption{\meshtag{Cartesian}.}
    \label{fig:mesh_a}
  \end{subfigure}%
  \begin{subfigure}[t]{0.33\textwidth}
    \centering \includegraphics[scale=0.25]{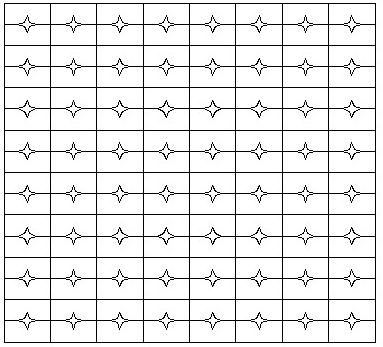}
    \caption{\meshtag{Star}.}
    \label{fig:mesh_b}
  \end{subfigure}
  \begin{subfigure}[t]{0.33\textwidth}
    \centering \includegraphics[scale=0.287]{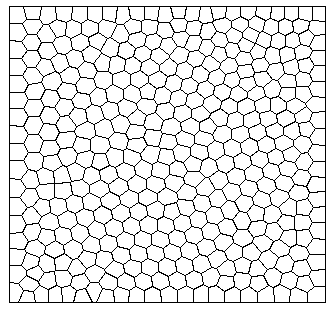}
    \caption{\meshtag{Polymesher}.}
    \label{fig:mesh_c}
  \end{subfigure}
  \caption{Test 1: meshes.}
  \label{fig:mesh}
\end{figure}

\begin{figure}
  \centering
  \begin{subfigure}[b]{0.33\textwidth}
    \centering \includegraphics[width=\linewidth]{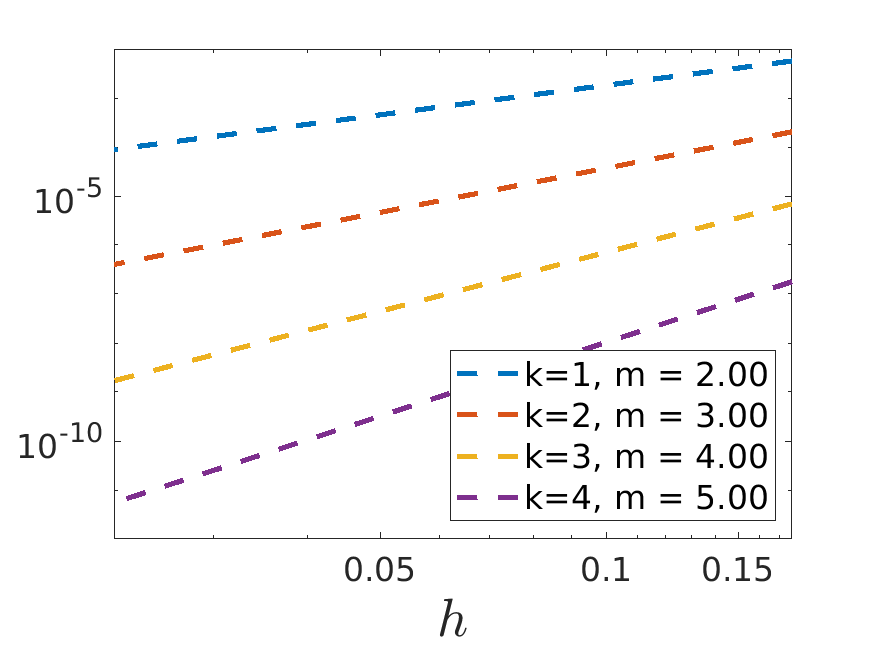}
    \caption{\meshtag{Cartesian}, $\lebl{}$ error.}
    \label{fig:state_a}
  \end{subfigure}%
  \begin{subfigure}[b]{0.33\textwidth}
    \centering \includegraphics[width=\linewidth]{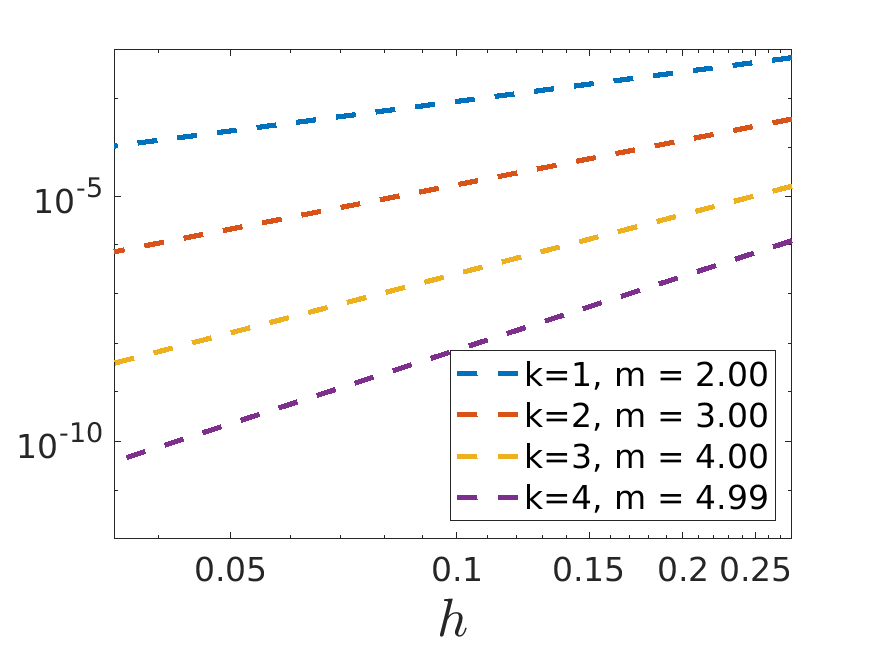}
    \caption{\meshtag{Star}, $\lebl{}$ error.}
    \label{fig:state_b}
  \end{subfigure}
  \begin{subfigure}[b]{0.33\textwidth}
    \centering \includegraphics[width=\linewidth]{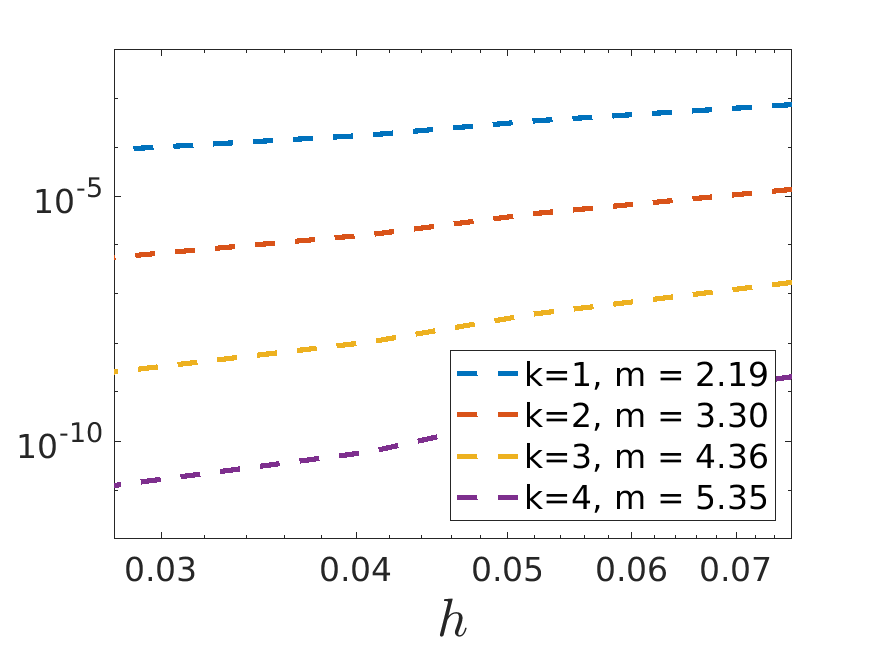}
    \caption{\meshtag{Polymesher}, $\lebl{}$ error.}
    \label{fig:state_c}
  \end{subfigure}
  \begin{subfigure}[b]{0.33\textwidth}
    \centering \includegraphics[width=\linewidth]{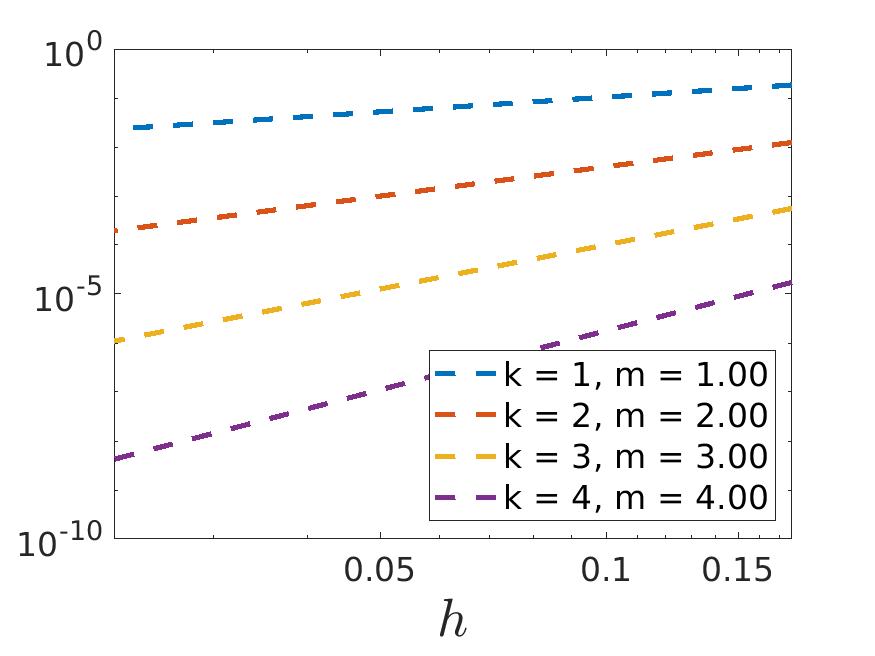}
    \caption{\meshtag{Cartesian}, e-error.}
    \label{fig:state_d}
  \end{subfigure}%
  \begin{subfigure}[b]{0.33\textwidth}
    \centering \includegraphics[width=\linewidth]{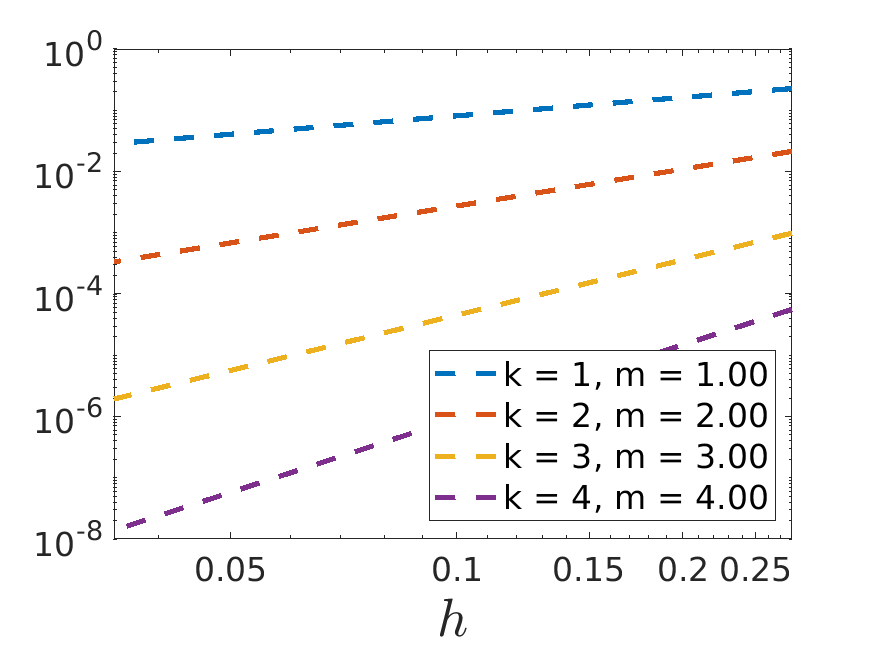}
    \caption{\meshtag{Star}, e-error.}
    \label{fig:state_e}
  \end{subfigure}
  \begin{subfigure}[b]{0.33\textwidth}
    \centering \includegraphics[width=\linewidth]{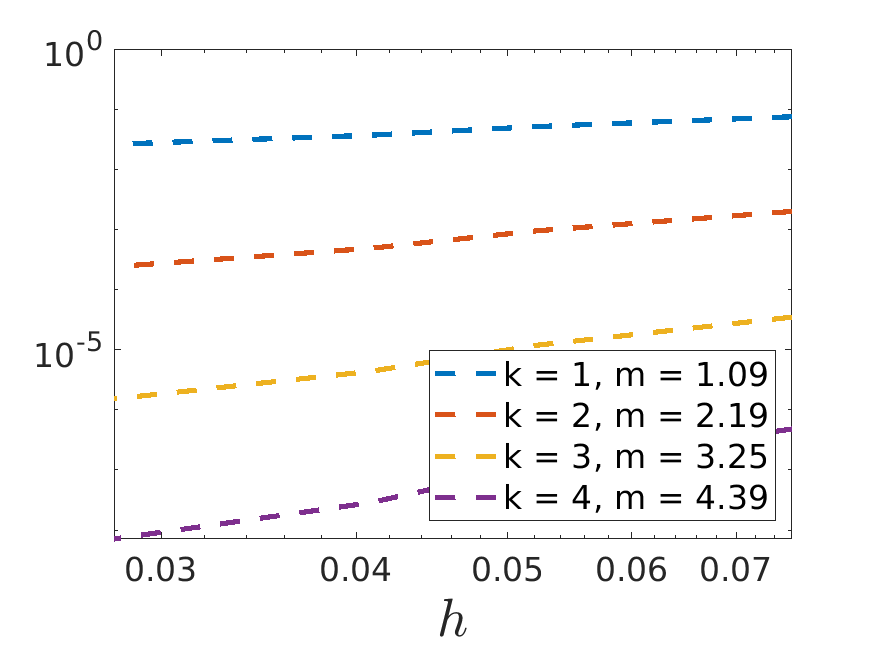}
    \caption{\meshtag{Polymesher}, e-error.}
    \label{fig:state_f}
  \end{subfigure}
  \caption{Test 1: convergence plots for the state variable $y$.}
  \label{fig:test1_state}
\end{figure}
\begin{figure}
  \centering
  \begin{subfigure}[b]{0.33\textwidth}
    \centering \includegraphics[width=\linewidth]{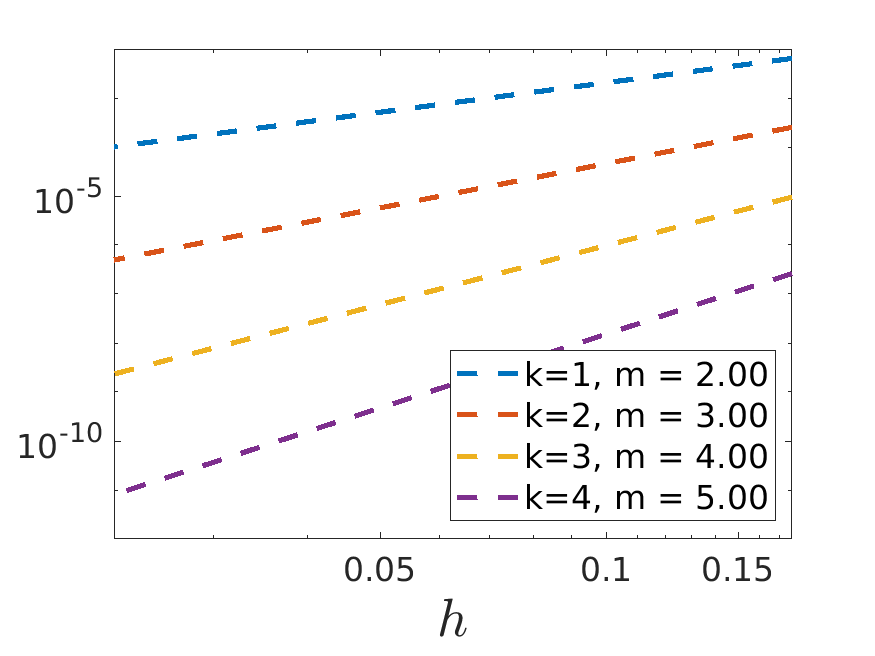}
    \caption{\meshtag{Cartesian}, $\lebl{}$ error.}
    \label{fig:adj_a}
  \end{subfigure}%
  \begin{subfigure}[b]{0.33\textwidth}
    \centering \includegraphics[width=\linewidth]{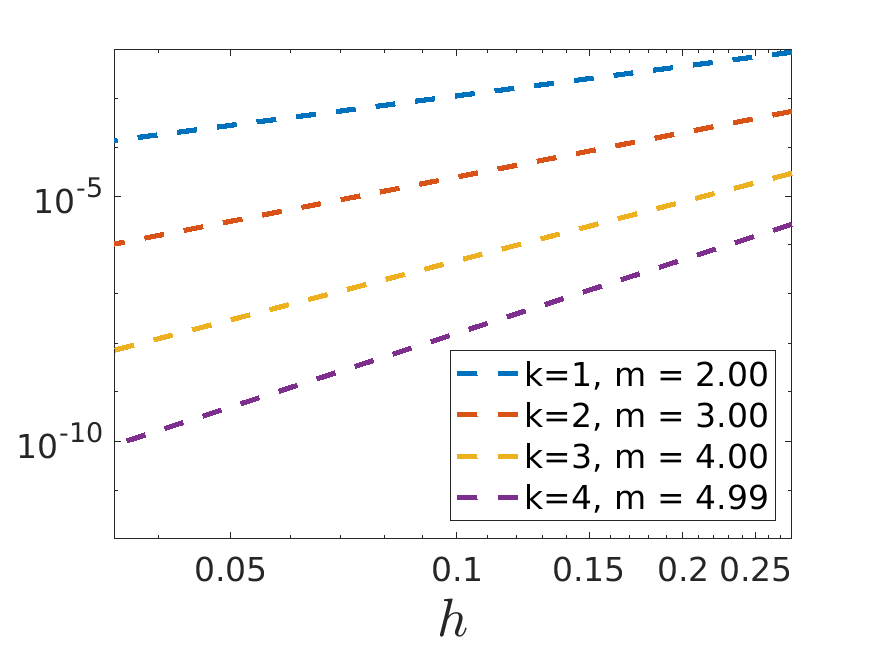}
    \caption{\meshtag{Star}, $\lebl{}$ error.}
    \label{fig:adj_b}
  \end{subfigure}
  \begin{subfigure}[b]{0.33\textwidth}
    \centering \includegraphics[width=\linewidth]{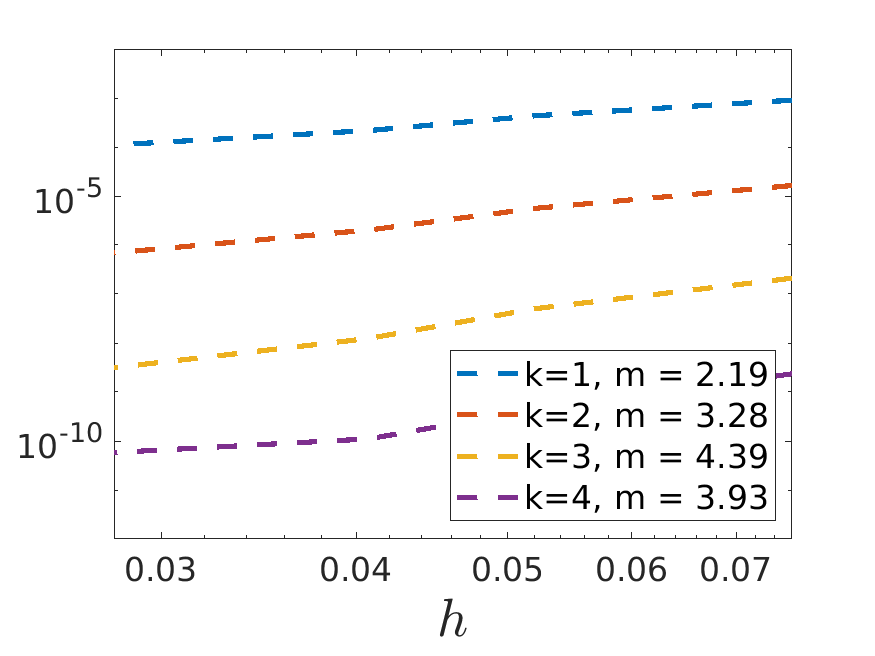}
    \caption{\meshtag{Polymesher}, $\lebl{}$ error.}
    \label{fig:adj_c}
  \end{subfigure}
  \begin{subfigure}[b]{0.33\textwidth}
    \centering \includegraphics[width=\linewidth]{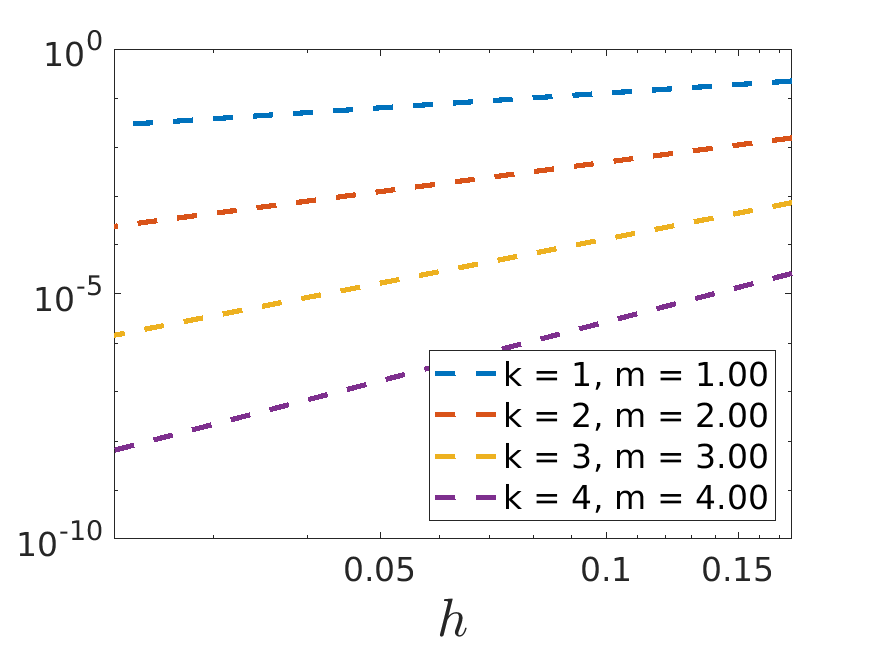}
    \caption{\meshtag{Cartesian}, e-error.}
    \label{fig:adj_d}
  \end{subfigure}%
  \begin{subfigure}[b]{0.33\textwidth}
    \centering \includegraphics[width=\linewidth]{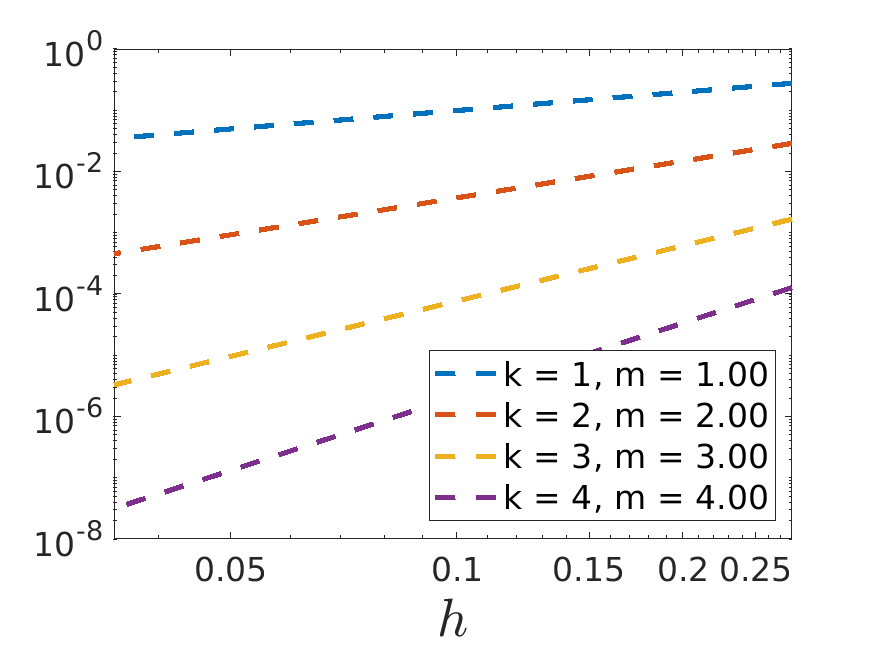}
    \caption{\meshtag{Star}, e-error.}
    \label{fig:adj_e}
  \end{subfigure}
  \begin{subfigure}[b]{0.33\textwidth}
    \centering \includegraphics[width=\linewidth]{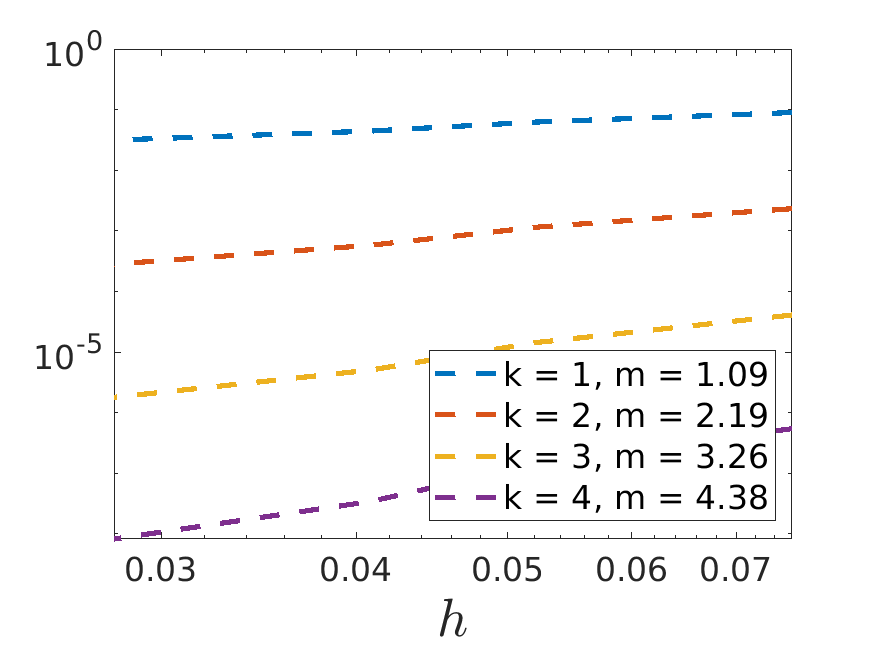}
    \caption{\meshtag{Polymesher}, e-error.}
    \label{fig:adj_f}
  \end{subfigure}
  \caption{Test 1: convergence plots for the adjoint variable $p$.}
  \label{fig:test1_adjoint}
\end{figure}
\begin{figure}
  \centering
  \begin{subfigure}[b]{0.33\textwidth}
    \centering \includegraphics[width=\linewidth]{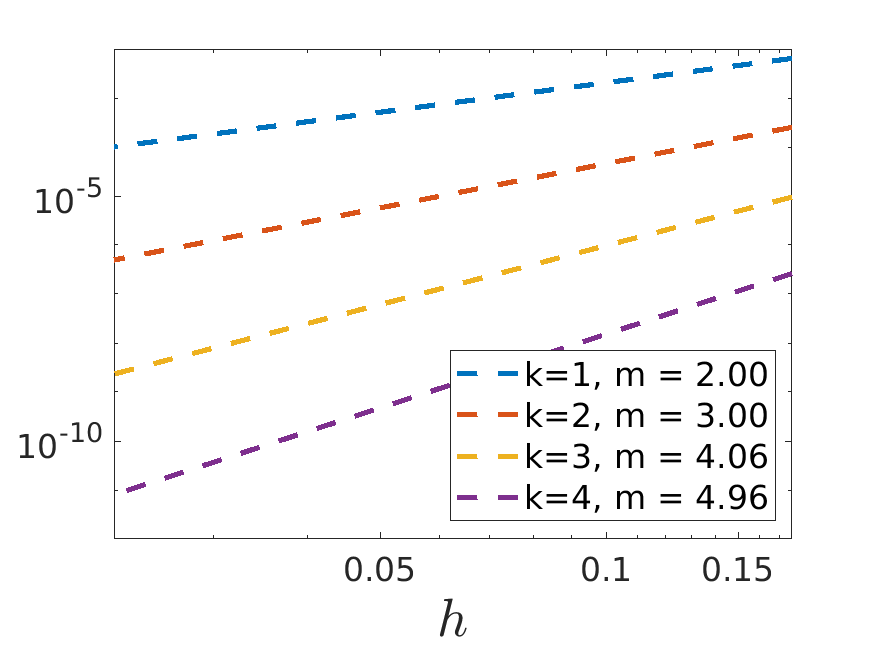}
    \caption{\meshtag{Cartesian}, $\lebl{}$ error.}
    \label{fig:control_a}
  \end{subfigure}%
  \begin{subfigure}[b]{0.33\textwidth}
    \centering
    \includegraphics[width=\linewidth]{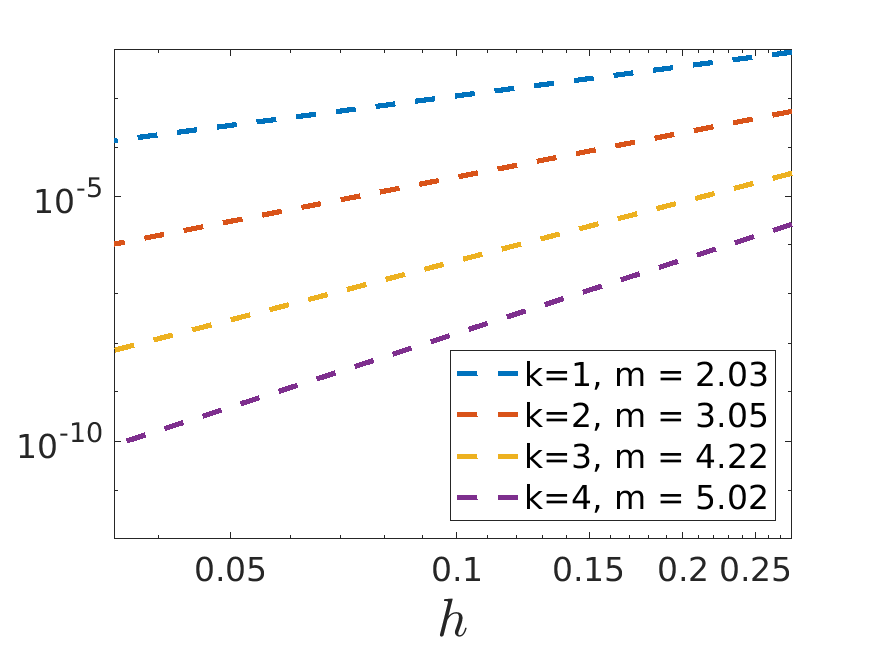}
    \caption{\meshtag{Star}, $\lebl{}$ error.}
    \label{fig:control_b}
  \end{subfigure}
  \begin{subfigure}[b]{0.33\textwidth}
    \centering \includegraphics[width=\linewidth]{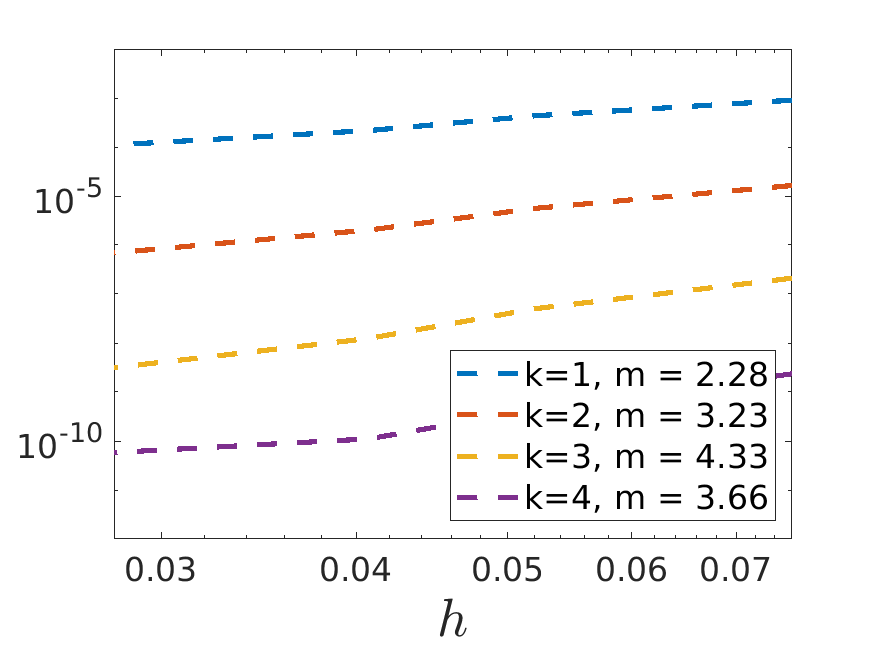}
    \caption{\meshtag{Polymesher}, $\lebl{}$ error.}
    \label{fig:control_c}
  \end{subfigure}
  \caption{Test 1: convergence plots for the control variable $u$.}
  \label{fig:test1_control}
\end{figure}

We solve Problem \ref{PB:VEM-SP} for $k=1,\ldots,4$ on three different grids,
depicted in Figure \ref{fig:mesh}, obtaining the approximated solutions
$y_h, u_h$ and $p_h$. 
The considered meshes are: a tessellation of
squares  (named \meshtag{Cartesian}, see Figure \ref{fig:mesh_a}), a mesh of non-convex octagons and
nonagons (named \meshtag{Star}, Figure \ref{fig:mesh_b}) and mesh of convex polygons obtained exploiting
\emph{Polymesher} \cite{Polymesher} (named \meshtag{Polymesher}, Figure \ref{fig:mesh_c}). 
For each mesh, we
consider four refinement levels to show the convergence properties stated in Theorem
\ref{thm:ErrorEstimate}, computing the following errors:
\begin{itemize}
\item the $\text{L}^2$ error measures
$$\displaystyle\sqrt{\sum_{E\in\Mh}\norm[\lebl{E}]{\proj{k}{E}y_h-y}^2}, \quad 
\displaystyle\sqrt{\sum_{E\in\Mh}\norm[\lebl{E}]{\proj{k}{E}p_h-p}^2}$$
and 
$$ \norm[\lebl{\Gamma_C}]{u_h - u},$$
\item the energy error (e-error) measures
$$\displaystyle\sqrt{ \sum_{E\in\Mh}\norm[\lebl{E}]{\sqrt{\K}
    \left ( \nabla \proj[\nabla]{k}{E}y_h - y \right
    )}^2
    },$$ and
$$ 
\displaystyle\sqrt{ \sum_{E\in\Mh}\norm[\lebl{E}]{\sqrt{\K} \left ( \nabla
      \proj[\nabla]{k}{E}p_h - p \right )}^2
    }.$$
\end{itemize}
We report the trend of the L$^2$ and the energy error for the state variable $y$ and the adjoint variable $p$
in Figures \ref{fig:test1_state} and \ref{fig:test1_adjoint}. The convergence plot for the L$^2$ error
of the control variable $u$ is depicted in Figure \ref{fig:test1_control}.
In the legends, we specify the average convergence rate
$m$ with respect to the mesh diameter $h$. 
We notice that the rates match the expected values stated in Theorem \ref{thm:ErrorEstimate}, even though some results for order $k=4$
are affected by machine precision issues.

\subsection{Test 2 (sensitivity study on the stabilization parameter $\sigma$)}
\begin{figure}
  \centering
  \begin{subfigure}[b]{0.33\textwidth}
    \centering \includegraphics[width=\linewidth]{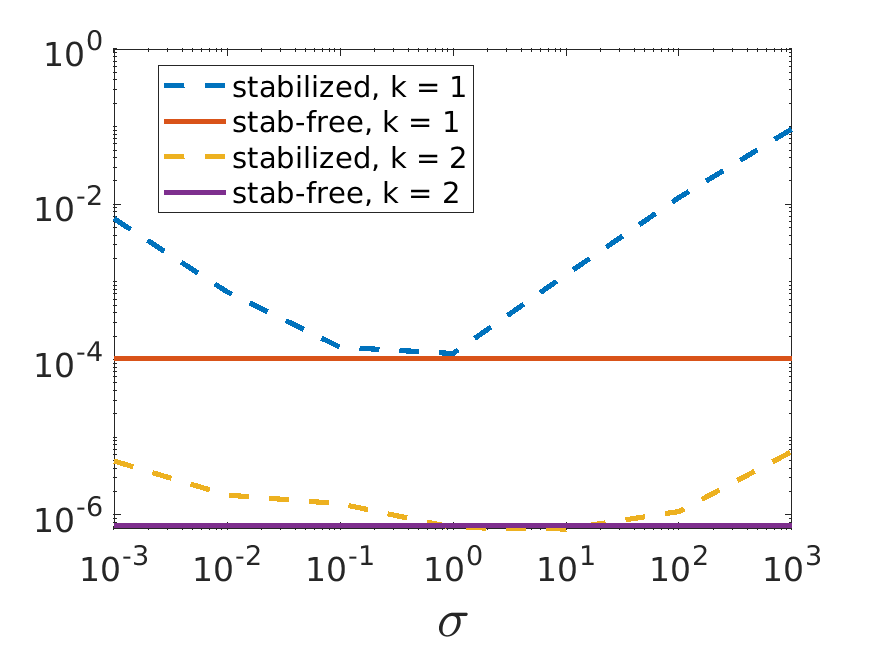}
    \caption{State, $L^2$ error.}
    \label{fig:test2_state_l2_starconcave}
  \end{subfigure}%
  \begin{subfigure}[b]{0.33\textwidth}
    \centering \includegraphics[width=\linewidth]{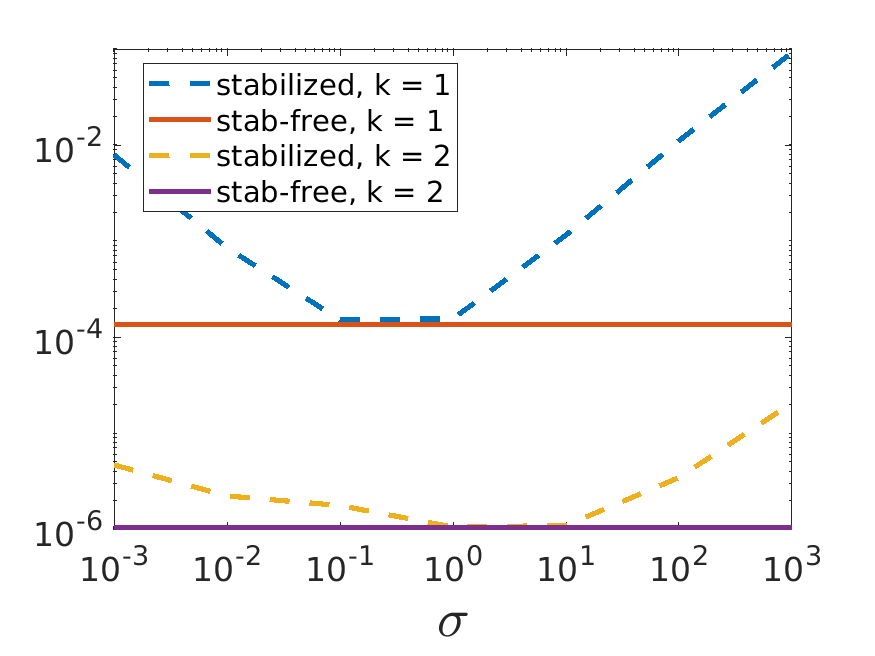}
    \caption{Adjoint, $L^2$ error.}
    \label{fig:test2_adj_l2_starconcave}
  \end{subfigure}
  \begin{subfigure}[b]{0.33\textwidth}
    \centering \includegraphics[width=\linewidth]{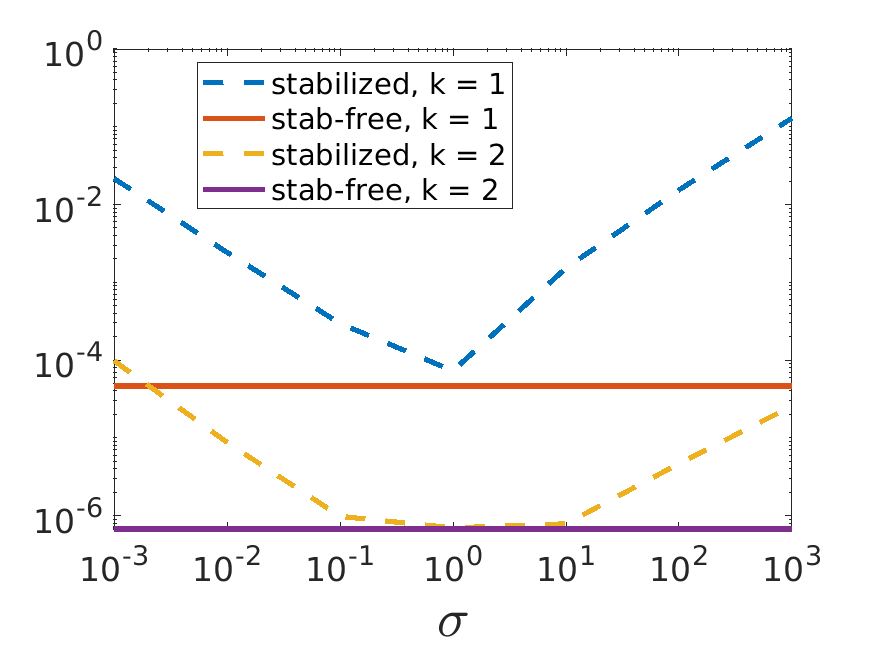}
    \caption{Control, $L^2$ error.}
    \label{fig:test2_control_l2_starconcave}
  \end{subfigure}
  \begin{subfigure}[b]{0.33\textwidth}
    \centering \includegraphics[width=\linewidth]{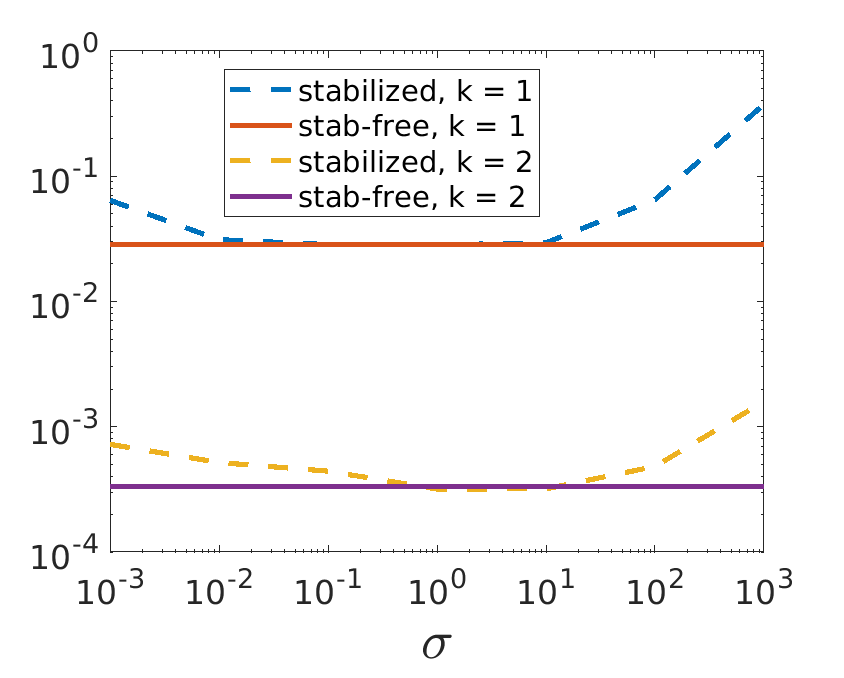}
    \caption{State, e-error.}
    \label{fig:test2_state_h1_starconcave}
  \end{subfigure}%
  \begin{subfigure}[b]{0.33\textwidth}
    \centering \includegraphics[width=\linewidth]{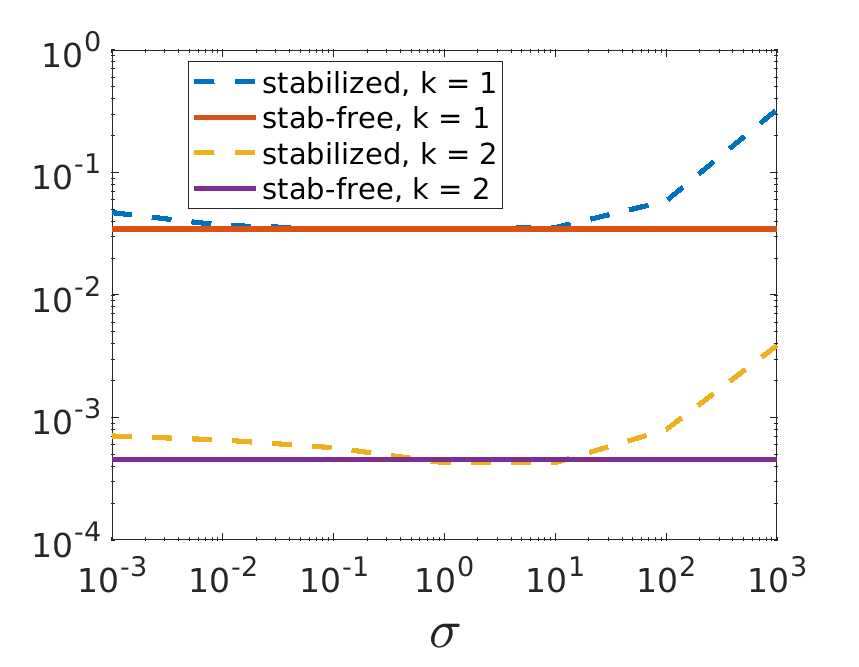}
    \caption{Adjoint, e-error.}
    \label{fig:test2_adj_h1_starconcave}
  \end{subfigure}
  \caption{\meshtag{Star}. Errors with respect to reference solutions, for $\sigma = 10^{-3},10^{-2}, 10^{-1}, 1, 10, 100, 1000$ and $k=1,2$.}
  \label{fig:test2_starconcave}
\end{figure}
\begin{figure}
  \centering
  \begin{subfigure}[b]{0.33\textwidth}
    \centering \includegraphics[width=\linewidth]{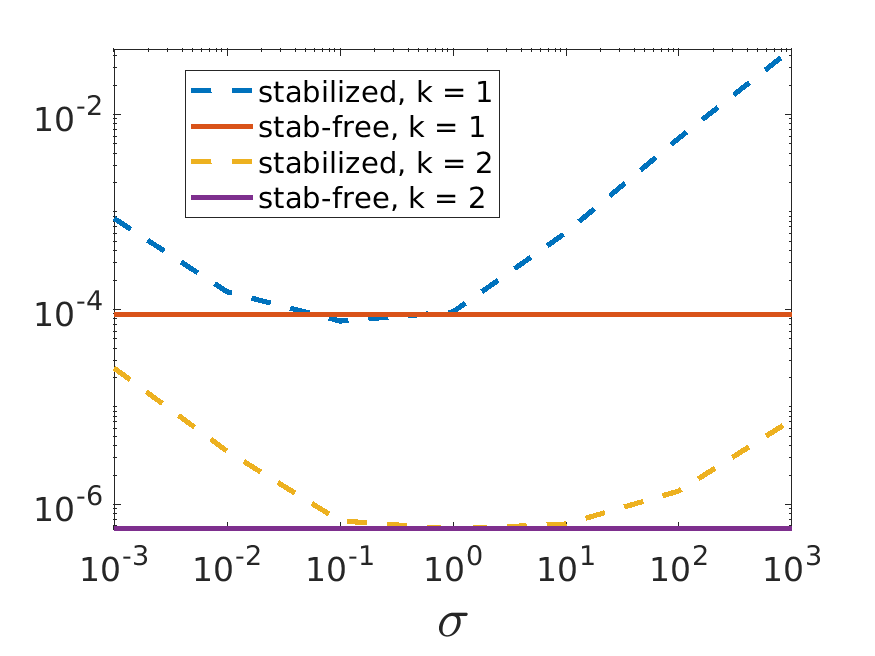}
    \caption{State, $L^2$ error}
    \label{fig:test2_state_l2_polymesher}
  \end{subfigure}%
  \begin{subfigure}[b]{0.33\textwidth}
    \centering \includegraphics[width=\linewidth]{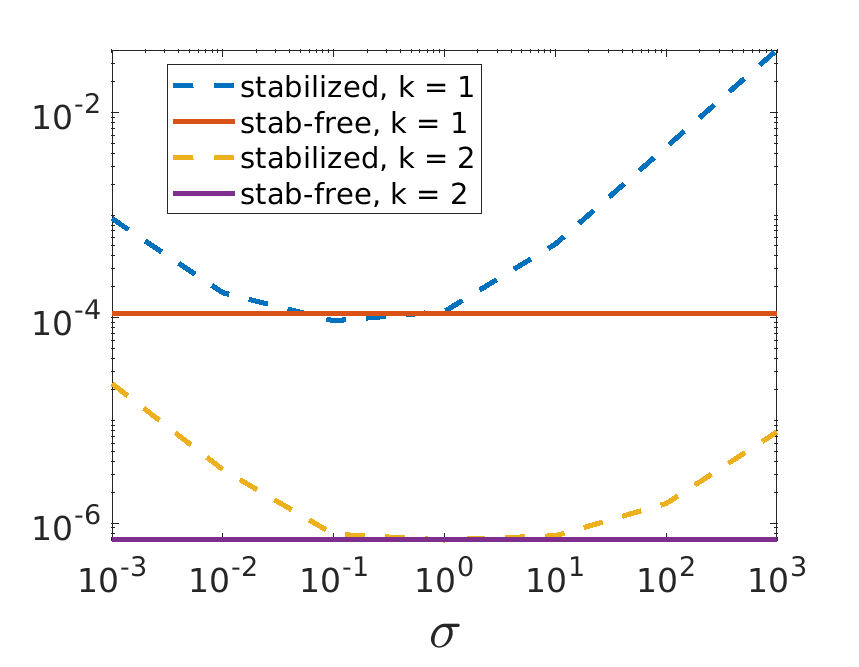}
    \caption{Adjoint, $L^2$ error}
    \label{fig:test2_adj_l2_polymesher}
  \end{subfigure}
  \begin{subfigure}[b]{0.33\textwidth}
    \centering \includegraphics[width=\linewidth]{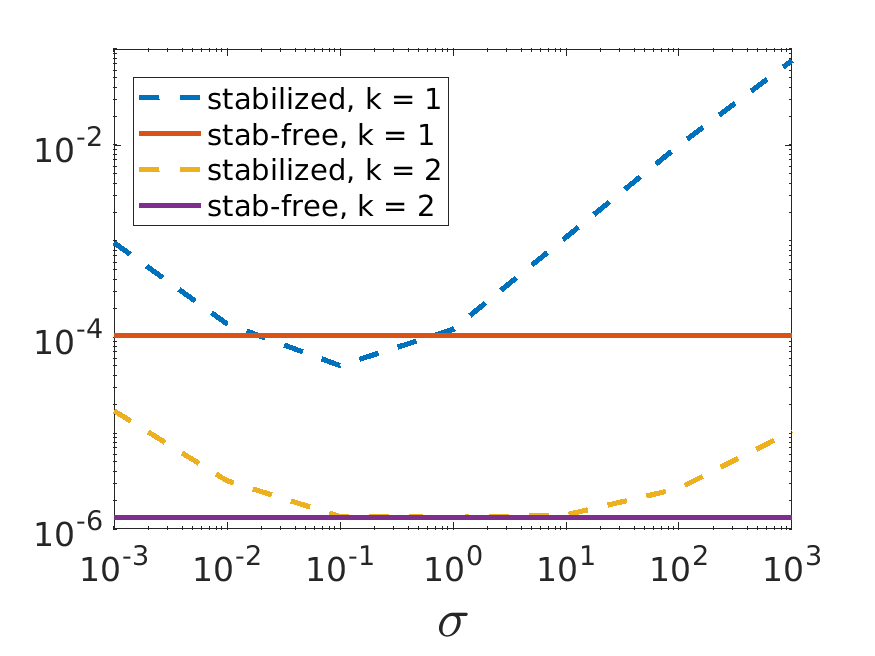}
    \caption{Control, $L^2$ error}
    \label{fig:test2_control_l2_polymesher}
  \end{subfigure}
  \begin{subfigure}[b]{0.33\textwidth}
    \centering 
    \includegraphics[width=1\linewidth]{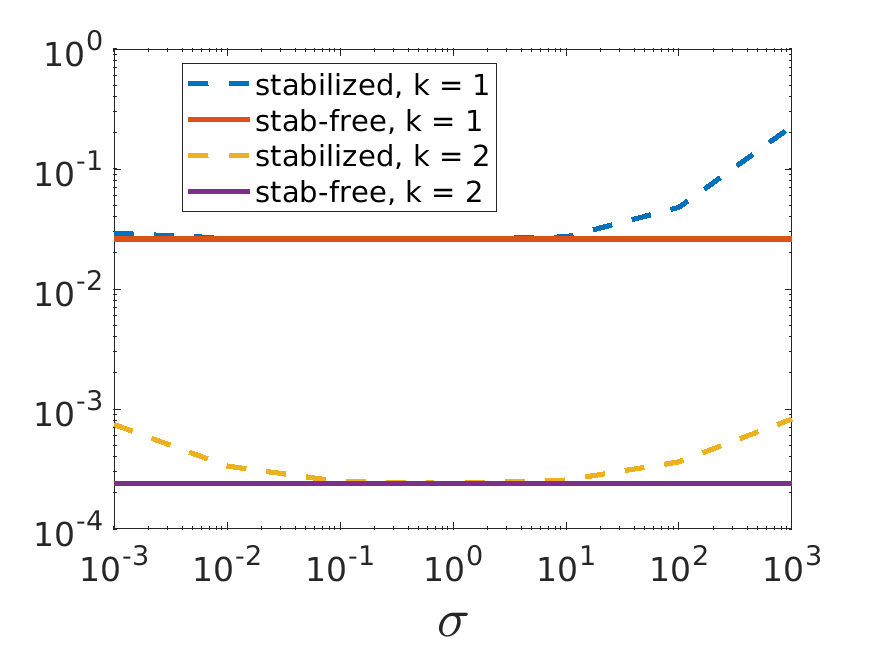}
    \caption{State, e-error}
    \label{fig:test2_state_h1_polymesher}
  \end{subfigure}%
  \begin{subfigure}[b]{0.33\textwidth}
    \centering \includegraphics[width=\linewidth]{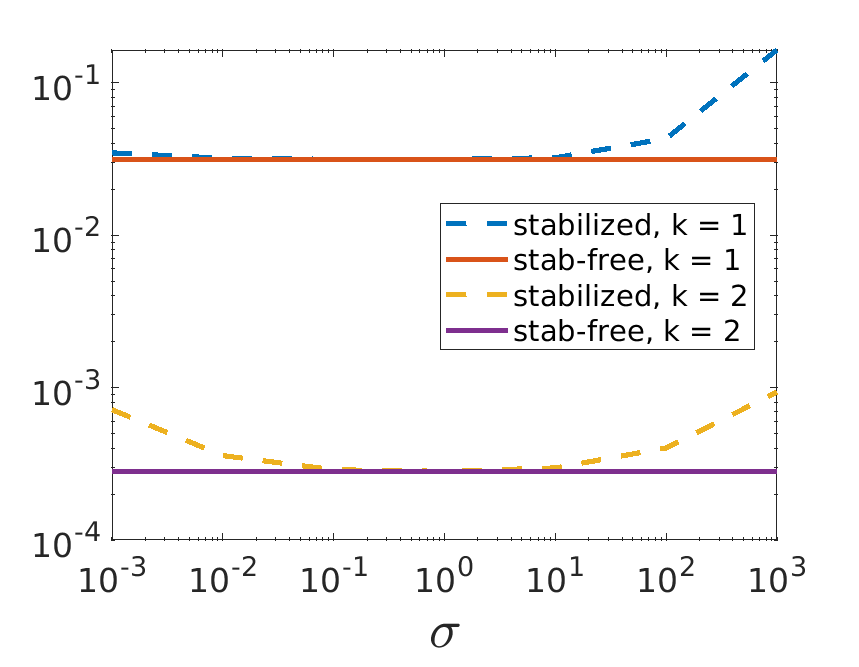}
    \caption{Adjoint, e-error}
    \label{fig:test2_adj_h1_polymesher}
  \end{subfigure}
  \caption{\meshtag{Polymesher}. Errors with respect to reference solutions, for $\sigma = 10^{-3},10^{-2}, 10^{-1}, 1, 10, 100, 1000$ and $k=1,2$.}
  \label{fig:test2_polymesher}
\end{figure}

In this section we consider the same problem as in Section \ref{sec:test1} and perform a numerical study of the errors done by the stabilized VEM scheme described in Remark \ref{rem:stab} for several choices of stabilization parameter $\sigma\in\{10^{-3},10^{-2}, 10^{1}, 1, 10, 100, 1000\}$. These results are compared with the results of the proposed stabilization-free scheme. In particular, we consider the most refined meshes belonging to the families \meshtag{Star} and \meshtag{Polymesher} depicted in Figure \ref{fig:mesh}, both featuring around 3000 polygons. In Figures \ref{fig:test2_starconcave} and \ref{fig:test2_polymesher} we display the $L^2$ and energy errors corresponding to the discrete solutions for $k=1,2$. We can observe a significant dependence of the error of stabilized VEM on the stabilization parameter. Moreover, the values of $\sigma$ corresponding to minima depend on the mesh (see in particular the results for the control variable). Finally, we observe that the errors made by the proposed stabilization-free scheme are always very close to the optimal errors made by stabilized VEM, independently of the solution and of the mesh. This suggests that the proposed method is a valid alternative to avoid issues related to the tuning of the stabilization parameter.

\subsection{Test 3 (towards applications)}
\label{sec:test2}
\begin{figure}
    \centering
    \begin{tikzpicture}[scale=3]
              \draw[color=black, thick, dashed] (2 +1.,0) --(2 + 2 ,0);
      \draw[color=black, thick, dashed] (2 +1.,1) --(2 + 2 ,1);
      \draw[color=black] (2 +1.,0.8) --(2 + 2 ,0.8); \draw[color=black] (2
      +1.,0.2) --(2 + 2 ,0.2); \draw[color=black, thick, densely dotted] (2
      +2.,1) --(2 + 2 ,0); \draw[color=black] (2 +1.,1) --(2 + 1 ,0.8);
      \draw[color=black] (2 +1.,.2) --(2 + 1 ,0.); \draw[color=black, very
      thick] (2 +1.,0) --(2 + 0 ,0) -- (2 + 0 ,1) -- (2 + 1 ,1);

      \node at (2 + .5, .5){\color{black}{$\Omega$}}; \node at (2 +
      1.5,0.9){$\Omega_1$}; \node at (2 + 1.5,0.1){$\Omega_2$}; \node at
      (2-.15,.5){\color{black}{$\Gamma_D$}}; \node at
      (2+2.15,.5){\color{black}{$\Gamma_N$}}; \node at (2 +
      1.2,1.1){$\Gamma_C$}; \node at (2 + 1.2,-.1){$\Gamma_C$};

      \node at (2 + 0,-.1){\color{black}{$(0,0)$}}; \node at (2 +
      2,-.1){\color{black}{$(2,0)$}}; \node at (2 +
      2,1.1){\color{black}{$(2,1)$}}; \node at (2 +
      0,1.1){\color{black}{$(0,2)$}};
    \end{tikzpicture}
    \caption{Test 3: schematic representation of the spatial domain. 
    }
    \label{fig:Domain_test2}
\end{figure}

\begin{figure}
  \centering
  \begin{subfigure}[t]{.49\linewidth}
    \includegraphics[width=\linewidth]{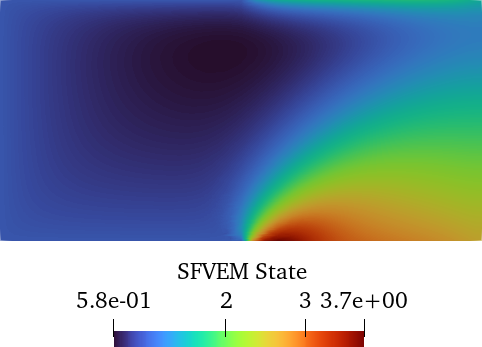}
  \end{subfigure}
  \hfill
    \begin{subfigure}[t]{.49\linewidth}
    \includegraphics[width=\linewidth]{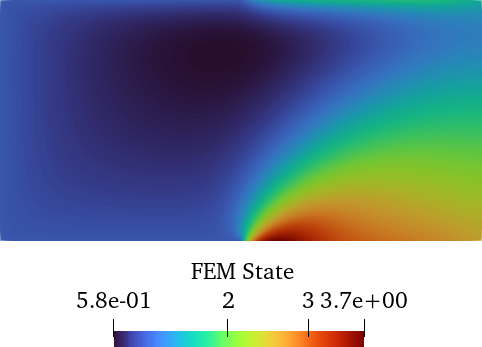}
  \end{subfigure}
  \\
  \vspace{0.1cm}
  \begin{subfigure}[t]{.49\linewidth}
    \includegraphics[width=\linewidth]{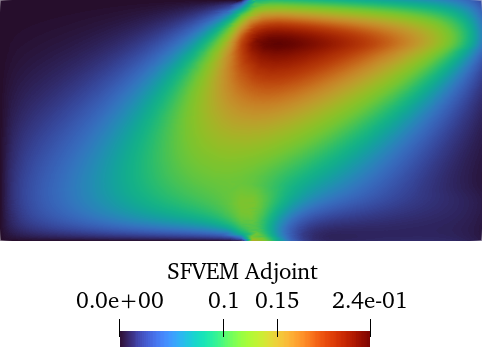}
  \end{subfigure}
  \hfill
    \begin{subfigure}[t]{.49\linewidth}
    \includegraphics[width=\linewidth]{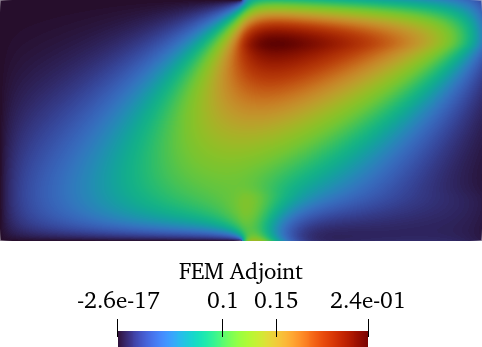}
  \end{subfigure}
    \caption{Test 3: Comparison of state and adjoint solutions obtained with the proposed method and the linear finite element implemented in FEniCS.}
  \label{fig:Test3StateAndAjoint}
\end{figure}

\begin{figure}
  \centering
  \begin{subfigure}[t]{.48\linewidth}
   \includegraphics[width=\linewidth]{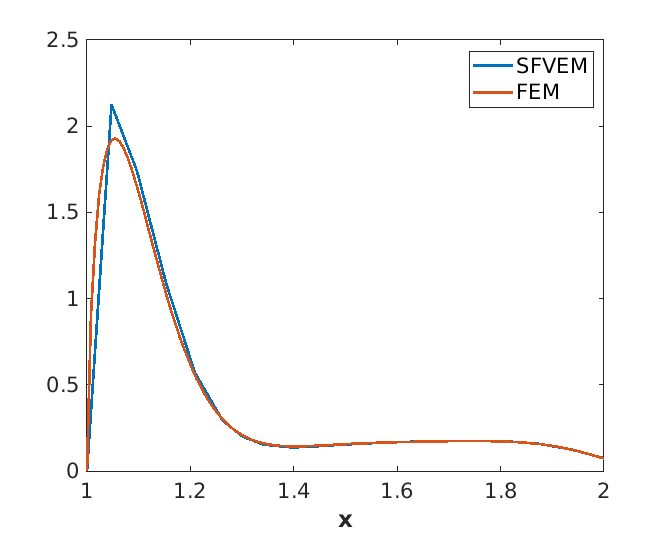}
    \caption{Control at $y=0$, $x\in[1,2]$.}
  \end{subfigure}
  \hfill
  \begin{subfigure}[t]{.48\linewidth}
   \includegraphics[width=\linewidth]{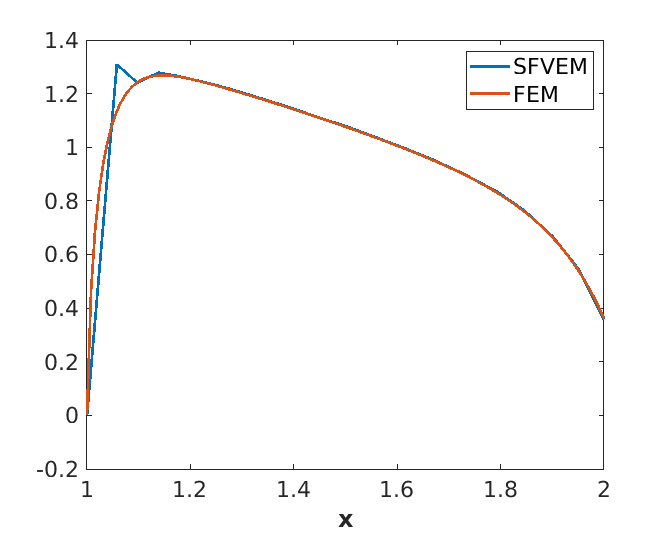}
    \caption{Control at $y=1$, $x\in[1,2]$.}
  \end{subfigure}
  \caption{Test 3: reference solution.}
  \label{fig:Test3Control}
\end{figure}

We propose a numerical study of the robustness of the proposed method, in its lowest-order formulation, by comparing it with a linear Finite Element discretization performed by the FEniCS library \cite{fenics} on the following OCP, in a setting inspired by \cite{negri2013reduced}.
The physical domain is $\Omega = (0,2) \times
(0,1)$. The desired state is observed in $\OmegaObs = \Omega_1 \cup \Omega_2$,
for $\Omega_1 = (1, 2)\times (0.8, 1)$ and $\Omega_2 = (1, 2)\times (0,
0.2)$. The control boundary $\Gamma_C$ is
$((1, 2) \times \{1\}) \cup ((1, 2) \times \{0\}).$ Homogeneous Neumann boundary
conditions are imposed over $\{2\}\times (0,2)$. On
$\Gamma_D = \partial \Omega \setminus (\Gamma_C \cup \Gamma_N)$, we impose
non-homogeneous Dirichlet boundary condition $g=1$. We remark that, even if the state
equation features non-homogeneous Dirichlet boundary conditions, the adjoint
equation still features homogeneous Dirichlet conditions, see \cite{manzoniOCP}.
The setting of this problem is represented in  Figure
\ref{fig:Domain_test2}.  
We solve Problem \ref{PB:VEM-SP} with $\alpha=0.07$,
$\mathcal K= 1/12$, $\gamma= 1$, $\beta=[1,1]$, $f=0$ and $\dState\equiv 2.5$ on a conforming (w.r.t $\OmegaObs$) mesh built using Polymesher \cite{Polymesher}, with $5490$ points.
We compare the SFVEM solution to the FEniCS reference FEM solution
over a triangular mesh with $28130$ points. 
In Figure \ref{fig:Test3StateAndAjoint}, we report the results of the state and the adjoint variables obtained by the SFVEM method and the FEM method. The SFVEM solution matches the FEM reference. 
Figure \ref{fig:Test3Control} shows the comparison of the solutions related to the control. The SFVEM control follows the behavior of the FEM reference, except for some small numerical oscillations around $x = 1$ due to the steep gradients in the solutions. {
This limitation can be effectively mitigated by refining the computational mesh in the region of interest. A more efficient strategy would involve an adaptive mesh refinement procedure, guided by an appropriate a posteriori error estimator. The development of this adaptive framework represent a promising direction for a future work.
}

\section{Conclusions}
\label{sec:conclusions}
In this work, we focused on a Stabilization-Free VEM approximation of Neumann boundary OCPs. We
addressed the problem in saddle point formulation, providing a rigorous \emph{a priori} error
estimation for all the variables at hand.  A numerical test based on analytic solutions confirmed
the expected theoretical convergence rates.  In a second test, we analyzed the role of the
stabilization term in the standard VEM formulation and its influence on the approximation error of
each component for this coupled problem, comparing it with the error made by the SFVEM scheme. For
standard VEM, we observed mesh-dependency effects and variabilities between the three variables of
the system, which are well circumvented by the SFVEM approach.  In the last test, we showed the
accuracy and robustness of the proposed method on a more application-oriented test.  This last test
is meant to be a bridge between the theoretical and applied aspects of the problem, where no exact
solution is known. Indeed, this work is a first step towards more complex applications.  There are
many aspects to be investigated to better show the SFVEM potential for such coupled problems, for
instance, more complicated geometrical structures and models, and the development of adaptive
strategies.

\section*{Acknowledgements}

The four authors are members of the Gruppo Nazionale Calcolo Scientifico (GNCS) at Istituto
Nazionale di Alta Matematica (INdAM). The authors kindly acknowledge financial support by INdAM-GNCS
Project 2025 (CUP: E53C24001950001).
AB and FM acknowledge funding by the Italian Ministry of University and Research (MUR)
through the project MUR-M4C2-1.1-PRIN 2022 
CUP: E53D23005820006. AB ackowledges funding by the European Union through Next Generation EU, M4C2, PRIN 2022 PNRR project CUP: E53D23017950001 and PNRR M4C2 project of CN00000013 National Centre for HPC, Big Data and Quantum
Computing (HPC) CUP: E13C22000990001.
\noindent MS thanks the ``20227K44ME - Full and Reduced order modelling of coupled systems: focus on non-matching methods and automatic learning (FaReX)'' project – funded by European Union – Next Generation EU  within the PRIN 2022 program (D.D. 104 - 02/02/2022 Ministero dell’Università e della Ricerca). This manuscript reflects only the authors’ views and opinions and the Ministry cannot be considered responsible for them.

\bibliographystyle{elsarticle-num} 
\bibliography{bib}






\end{document}